\documentclass[12pt, reqno]{amsart}
\usepackage[utf8]{inputenc}
\usepackage{amsmath,amsthm}
\usepackage{amssymb}
\usepackage{amsfonts}

\usepackage[left=2cm,right=2cm,top=3cm,bottom=2cm]{geometry}

\newtheorem{theorem}{Theorem}[section]
\newtheorem{lemma}[theorem]{Lemma}
\newtheorem{definition}[theorem]{Definition}
\newtheorem{corollary}[theorem]{Corollary}
\newtheorem{proposition}[theorem]{Proposition}
\newtheorem{claim}[theorem]{Claim}

\theoremstyle{definition}
\newtheorem{example}[theorem]{Example}

\newcommand\blfootnote[1]{
  \begingroup
  \renewcommand\thefootnote{}\footnote{#1}
  \addtocounter{footnote}{-1}
  \endgroup
}

\begin{document}

\renewcommand{\refname}{References}
\thispagestyle{empty}

\begin{center}{\bf
CONSTRUCTING MODELS OF SMALL ORDERED THEORIES WITH MAXIMAL COUNTABLE SPECTRUM}
\end{center}

\noindent{}

\begin{center}
BEKTUR BAIZHANOV AND TATYANA ZAMBARNAYA
\blfootnote{Second author was supported by the 
Committee of Science of the Ministry of Education and Science of the Republic of Kazakhstan, Grant No. AP08955727.}
\end{center}

\noindent{}

\begin{center}
\begin{minipage}{0.85\textwidth}
\begin{footnotesize}
\textsc{Abstract.}
We present a method for constructing countable models of small theories and
apply it to prove theorems on the maximal number of countable non-isomorphic models
of linearly ordered theories.
\end{footnotesize}
\end{minipage}
\end{center}

\section*{Introduction}

The Vaught's conjecture was confirmed for 
theories of pure linear order with a finite or a countable set of unary predicates by M.~Rubin~\cite{R1974},
for o-minimal theories by L.~Mayer~\cite{LM1988}, and for 
quite o-minimal theories by Sudoplatov-Kulpeshov~\cite{KulSud2017};
as well as for cases with the following restrictions:
for weakly o-minimal theories of convexity rank one by Alibek-Baizhanov-Kulpeshov-Zambarnaya \cite{ABKZ2018}, 
for binary, stationarily ordered theories by Moconja-Tanovic \cite{MC2019}, 
and for weakly o-minimal theories of finite convexity rank by B.~Kulpeshov \cite{KBSH2020}.

The first step in describing countable spectrum is determining when theories have 
the maximal number of countable non-isomorphic models.
In~\cite{Kud1986, Sud2004, Sud2018} K.Zh. Kudaibergenov and  S.V.~Sudoplatov 
proved that every countable model of a small theory can be represented
as a union of an elementary chain of prime models over finite sets.
Based on ideas from this theorem, in~\cite{ABZ2014, BBZ2018, BUZ2019a, BUZ2019b},
authors used Tarski-Vaught test to construct countable pairwise non-isomorphic models.
In section~\ref{Sconstruction} we present a generalization of this construction (Theorem~\ref{Construction}),
such that the obtained model contains a given countable set and realizes the smallest number of types.
Unlike Kudaibergenov-Sudoplatov's result, we build a countable elementary submodel of an $\aleph_1$-saturated model
satisfying the desired properties, rather than reconstruct a given structure.
The method results in an elementary chain of prime models over finite sets
and preserves non-homogeneity (Corollary~\ref{homogeneous1}).
The main idea of the construction is in choosing elements as realizations of decreasing
sequences of principal formulas in an $\aleph_1$-saturated model.
Later, in sections~\ref{Squasi}, \ref{Strivial} and \ref{SEG},
this construction is applied to simplify proofs of some theorems (\ref{abzQuasi}, \ref{Trivial} and \ref{NoGreatest})
 and prove new theorems (\ref{NoLeast}, \ref{Infinite}) and corollaries
on maximal number of countable non-isomorphic models  of small ordered theories
with a natural restriction (conditions) on 1-types and 2-formulas definable over finite sets.

\section{Preliminaries}

Further in the article we consider small theories,
that is, theories $T$ with $|\underset{n<\omega}{\bigcup}S_n(T)| = \aleph_0$.
By Gothic letters ($\mathfrak A$, $\mathfrak M$, $\mathfrak N$, ...) we denote structures,
and universes of those structures we denote by capital letters ($A$, $M$, $N$, ...).

For subsets $A$ and $ B$ of an ordered structure $\mathfrak M$  we use the following notations:

\begin{center}
$A^+:=\{\gamma\in M\ |$ for all $ a\in A,\ \mathfrak M\models a<\gamma\}$;

$A^-:=\{\gamma\in M\ |$ for all $ a\in A,\ \mathfrak M\models \gamma<a\}$.
\end{center}
We write $A<B$ if for all $a\in A$, $b\in B$ $\mathfrak M\models a<b$.
If $A$ and $B$ are $C$-definable ($C\subseteq M$), then $A^+, A^-$ and $A<B$ are $C$-definable as well.

For a 2-formula $\varphi(x,y)$ denote
$\varphi (x,y)^-:= \forall z \big(\varphi(z, y)\to z>x\big)$, and
$\varphi (x,y)^+:= \forall z \big(\varphi(z, y)\to z<x\big)$.

\begin{definition}
The set $A$ is said to be {\bf convex} in a set $B\supseteq A$ if for all
$a,b\in A$ and all $c\in B$ $a<c<b$ implies that $c\in A$.
The set $A$ is {\bf convex} if it is convex in $M$.
\end{definition}

\begin{definition}
A formula $\varphi(x,\bar y, \bar a)$ is a {\bf convex} formula, if for every $\bar b\in M$
the set $\varphi(M,\bar b,\bar a)$ is convex in every model of $Th(\mathfrak M)$ containing $\bar b$ and $\bar a$.
\end{definition}

\begin{definition}\cite{BV2011}
1) A {\bf convex closure} of a formula $\varphi(x,\bar a)$
is the following formula:

\begin{center}
$\varphi^c(x,\bar a):=\exists y_1\exists y_2 \big(\varphi(y_1,\bar a)\land\varphi(y_2,\bar a) \land (y_1\leq x\leq y_2)\big)$.
\end{center}

2) A {\bf convex closure} of a type $p(x)\in S_1(A)$ is the type

\begin{center}
$p^c(x):=\{\varphi^c(x,\bar a)\ |\ \varphi(x,\bar a)\in p\}$.
\end{center}
\end{definition}
Similarly, denote $tp^c (\alpha/A):=\{\varphi^c(x,\bar a)\ |\ \varphi(x,\bar a)\in tp (\alpha/A)\}$.

In weakly o-minimal and, therefore, o-minimal theories $p^c(\mathfrak M)=p(\mathfrak M)$ for every $p\in S_1(A)$.

\section{Construction of a countable model}\label{Sconstruction}

\begin{definition}
Let $\mathfrak M$ be a model of a theory $T$. A {\bf finite diagram}~\cite{Sh1970} of
$\mathfrak M$ is the collection of all $\emptyset$-definable complete types that are realized in $\mathfrak M$:
\begin{center}
$\mathcal D(\mathfrak M)=\{p\in S(T)\ |\ \mathfrak M\models p\}$.
\end{center}

A {\bf dowry} of a set $B\subseteq M$ is the following collection:

\begin{center}
$\mathcal D(B)=\{p\in S(T)\ |$ there exist $\mathfrak M\models T$ and $\bar b\in B$ such that $\bar b\in p(\mathfrak M)\}$.
\end{center}

\end{definition}

We present a lemma from \cite{BBZ2018} (Lemma 3.5) together with its proof since later we refer to the construction in the proof.
The set $A$ from \cite{BBZ2018} is considered to be empty.

\begin{lemma}\label{chain}\cite{BBZ2018}
Let $T$ be a small countable complete theory,
$D$ be a finite and $B=\{b_1, b_2, ..., b_i,...\}$ be a countable subset of a model of $T$.
Then for each satisfiable $(B\cup D)$-formula $\varphi(x,\bar b_n,\bar d)$,
where $\bar b_n:=( b_1, b_2, ..., b_n)\in B$ ($n<\omega$) and $\bar d\in D$,
there exists a type $q_\varphi \in S_1(B\cup D)$ such that

1) $\varphi(x, \bar b_n,\bar d)\in q_\varphi$;

2) for every $i\geq n$ the type $q_\varphi\restriction (B_i\cup D)$ is principal,
where $B_i=\{b_1, b_2, ..., b_i\}$.
\end{lemma}

\noindent{}{\it Proof of Lemma~\ref{chain}.}
Let $\bar d'$ be a tuple enumerating the set $D$.
Because the theory $T$ is small, there exists a formula
$\varphi'(x, \bar b_n,\bar d)$ that implies $\varphi(x,\bar b_n,\bar d)$
and generates a principal type over $\{\bar b_n, \bar d\}$.
In turn there is a principal subformula $\varphi_0(x, \bar b_n,\bar d')$ over $(\{\bar b_{n}\}\cup D)$ that implies
$\varphi'(x, \bar b_n,\bar d)$, and, analogically,
a principal subformula over $(\{\bar b_{n+1}\}\cup D)$ that implies
$\varphi_0(x, \bar b_n,\bar d')$.
Repeating this procedure, we obtain a
consistent infinite decreasing chain of principal over parameters formulas
$\varphi_i(x,\bar b_{n+i},\bar d')$:
... $\subseteq\varphi_{i+1}(M,\bar b_{n+i+1},\bar d')\subseteq$
$\varphi_i(M,\bar b_{n+i},\bar d')\subseteq$ ...
$\subseteq\varphi_0(M,\bar b_n,\bar d')\subseteq$
$\varphi'(M,\bar b_n,\bar d)\subseteq\varphi(M,\bar b_n,\bar d)$,
where $\mathfrak M$ is an arbitrary model of $T$ with $(B\cup D)\subseteq M$.
By this we have defined the desired complete type over $(B\cup D)$.

\hfill $\Box$

\begin{theorem}[Tarski---Vaught Test]\cite{TZ2012}\label{TV}
Let $\mathfrak B$ be an $L$-structure and $A$ be a subset of $B$.
Then $A$ is the universe of an elementary substructure of $\mathfrak B$ if and only if every
$L(A)$-formula $\varphi(x)$ which is satisfiable in $\mathfrak B$ can be satisfied by an
element of $A$.
\end{theorem}

The following theorem gives us a method of construction of a countable model by use of
the Tarski---Vaught Test.

\begin{theorem}\label{Construction}
Let $\mathfrak M\models T$ be a model of a small countable complete theory $T$.
Let $B\subseteq M$ be countable.
Then there exists a countable model $\mathfrak A^B=\mathfrak A\models T$ such that $B\subseteq A$ and

1) for every $\bar a\in A$ there is $\bar b\in B$, such that $tp(\bar a/\bar b)$ is principal;

2) for every $\mathfrak C\models T$ with $B\subseteq C$, $\mathcal D(\mathfrak A)\subseteq\mathcal D(\mathfrak C)$;

3) all such models $\mathfrak A$ obtained by the given construction have the same finite diagram.

\end{theorem}

\noindent{}{\it Proof of Theorem~\ref{Construction}.}
Denote by $\mathfrak N$ an $\aleph_1$-saturated elementary extension of $\mathfrak M$.
Fix some enumeration $B=\{b_1, b_2, ..., b_i,...\}$.
We will use Tarski-Vaught criterion in order to show that the constructed set $A$ is the universe of an
elementary substructure of $\mathfrak N$.
On each step of the construction we will be fixing
a set of parameters and promising to realize each satisfiable 1-formula over it.
We must keep coming back to the same set of parameters
and deal with another formula. So, the different sets of parameters are being attacked in parallel.

{\bf Step 1.}
Denote by
$\Phi_1$
the set of all $\emptyset$-definable unary formulas,
$\Phi_1:=\{\varphi_i^1(x)\ |\ i<\omega\}$.
Choose the least $i$ such that
$\mathfrak N\models\exists x \varphi_i^1(x)$.
To satisfy the Tarski-Vaught property, we must find a witness for $\varphi_i^1(x)$.
If there exists a principal over $\emptyset$ subformula of $\varphi_i^1(x)$ that has a non-empty intersection with $B$,
choose the smallest index $j$ such that $b_j\in \varphi^1_{i,0}(N)\cap B$ for some principal subformula
$\varphi^1_{i,0}$ of $\varphi^1_{i}$.
Denote $d_1:=b_j$, $D_1:=\{d_1\}$, and $B_1:=\{b_1,d_1\}$.

Otherwise, since the set $B$ and the formula $\varphi_i^1$ are as in Lemma~\ref{chain}
(let the set $D$ to be empty),
there exists a $B$-type $q_{\varphi_i^1}$ satisfying conditions 1) and 2) of the lemma.
And since the model $\mathfrak N$ is $\aleph_1$-saturated, this type is realized in $\mathfrak N$ by
some element, denote it by $d_1$. Thus, $d_1$  is principal over $\emptyset$.
Denote $D_1:=\{d_1\}$, and $B_1:=\{b_1\}$.

{\bf Step 2.}
Choose the least $j$ such that the formula  $\varphi_j^1(x) \in \Phi_1$
was not considered before with
$\mathfrak N\models\exists x \varphi_j^1(x)$.
We find a special witness for $\varphi_j^1(x)$.
If exists, denote by $d_2$ an element with the smallest index in the enumeration of $B$
such that for some principal subformula $\varphi_j^1(x)$ of $\varphi_j^1(x)$
$d_2\in\varphi_j^1(N)\cap B$.
Denote $B_2:=B_1\cup\{b_2, d_2\}\subset B$.
Otherwise, apply Lemma~\ref{chain} to the sets $B$ and $\{d_1\}$, and the formula $\varphi_j^1(x)$,
to find a realization $d_2$ of the type $q_{\varphi_j^1}$, which exists by the lemma.
We can arrange that $d_2$ is principal over $d_1$.
Denote $B_2:=B_1\cup\{b_2\}$.
Then, $b_2$ will be principal over every superset of $B_2$.

Now, consider the set of all $B_1\cup D_1$-definable 1-formulas
$\Phi_2:=\{\varphi_i^2(x, b_1, d_1)\ |\ i<\omega\}$.
Choose  the least index $i$ such that the formula $\varphi_i^2(x, b_1, d_1)\in\Phi_2$
was not considered previously, and
$\mathfrak N\models\exists x \varphi_i^2(x, b_1, d_1)$, and find
a realization $d_3$ either as an element with the smallest index in $B$ that it is in
$\varphi_{i,0}^2(N, b_1, d_1)\cap B$ for some principal over $B_1\cup D_1$
 subformula $\varphi_{i,0}^2$ of $\varphi_{i}^2$, or,
by applying Lemma~\ref{chain} to
$B$, $\{d_1, d_2\}$ and $\varphi_i^2$.
In the first case, rename $B_2$ to be $B_2\cup\{d_3\}$.
Denote $D_2:=\{d_1, d_2, d_3\}$.

By the end of step $k$ the following sets will be defined:

1) Finite nested sets
$D_1=\{d_1\}$,   $D_2=\{d_1, d_2, d_3\}$,
$D_3=\{d_1, d_2, ..., d_6\}$, ...,   $D_k=\{d_1, d_2,...,$ $d_{\frac{(k+1)k}{2}}\}$, where
$D_i$ is obtained on step $i$ by adding $i$ realizations to the set $D_{i-1}$.
It is possible that $d_i=d_j$ for some $i$ and $j$ with $1\leq i<j\leq \frac{(k+1)k}{2}$.

2) Finite sets of parameters $B_1\subseteq B_2\subseteq...\subseteq B_{k}\subset B$,
where each $B_i$ is obtained from $B_{i-1}$ by adding $b_i$ and, maybe, some realizations of
formulas from $\Phi_1,...,\Phi_i$.

3) The family of all $\emptyset$-definable 1-formulas $\Phi_1$,
and, for every $i$, $2\leq i\leq k$,
the family of $B_{i-1}\cup D_{i-1}$-definable 1-formulas, $\Phi_i$.

Further we use the usual notation $\bar d_i=( d_1, d_2, ..., d_i )$, $i<\omega$.

{\bf Step $k+1$.}
Denote $B_{k+1}:=B_k\cup\{b_k+1\}$.
Firstly, we realize one formula from each of the families we defined earlier.
To do this, for each $m$, $1\leq m\leq k$, find smallest index $i_m$ such that
the  formula $\varphi_{i_m}^m\in \Phi_m$ was not considered before,
and definable set of which in the model $\mathfrak N$ is not empty.
If exists, choose the smallest index $j$ such that $b_j\in \varphi_{i_m,0}^m\cap B$ for
some principal over $B_{m-1}\cup D_{m-1}$ subformula $\varphi_{i_m,0}^m$ of $\varphi_{i_m}^m$.
Denote $d_{\frac{(k+1)k}{2}+m}:=b_j$, and rename $B_{k+1}$ to be $B_{k+1}\cup\{d_{\frac{(k+1)k}{2}+m}\}$.
Otherwise, apply Lemma~\ref{chain} to the sets $B$ and $\{d_1, d_2, ..., d_{\frac{(k+1)k}{2}+m-1}\}$, and
the formula $\varphi_{i_m}^m$, to find a realization $d_{\frac{(k+1)k}{2}+m}$ of the type
$q_{\varphi_{i_m}^m}$.
We can arrange that $d_{\frac{(k+1)k}{2}+m}$ is principal over $B_{m-1}$ and the $d_j$'s for $j<\frac{(k+1)k}{2}+m$.

Now, let $\Phi_{k+1}$ be the set of all $B_k\cup D_{k}$-definable
1-formulas, and find the smallest $i$ such that
$\mathfrak N\models \exists x\varphi^{k+1}_i(x, \bar b_k,\bar d_{\frac{(k+1)k}{2}})$
and that the formula $\varphi^{k+1}_i$ was not considered before.
We choose $d_{\frac{(k+1)k}{2}+k+1}$ as before, either as
an element of $B$ with the smallest index, that is in realization of some principal over $B_k\cup D_{k}$
subformula of $\varphi^{k+1}_i$,
or, as a realization of a type $q_{\varphi^{k+1}_i}$,
which exists by Lemma~\ref{chain} applied to
the sets $B$, $\{d_1, d_2, d_{\frac{(k+1)k}{2}+k}\}$ and formula $\varphi^{k+1}_i$.
In the first case, add $d_{\frac{(k+1)k}{2}+k+1}\in B$ to $B_{k+1}$.
Let $D_{k+1}$ be the set $\{d_1, d_2, ..., d_{\frac{(k+1)k}{2}+k+1}\}$.

Denote $A:= B\cup \underset{i<\omega}{\bigcup}D_i$.

The Tarski-Vaught criterion implies that the obtained set $A$
is a universe of an ele\-men\-tary substructure $\mathfrak A$ of $\mathfrak N$.

For $i<\omega$, when choosing $d_{i}$ to be as in Lemma~\ref{chain},
the type $tp(d_{i}/ B_n\cup\{\bar d_{i-1}\}$ is principal for every $n\geq i-1$.
If we choose $d_i\in B$ on step $j$ as a realization of a principal subformula, we add it to the set
$B_{j+1}$, as well as to all the next $B_{j+l}$'s.
This way, $d_i$ is principal over every superset of $B_{j+1}$.
Note that $d_1$ is principal over the empty set.
\begin{equation}\label{eqS}
\tag{*}
\begin{minipage}{0.9\textwidth}
From the previous statements it easily follows by induction that for every $i<\omega$ and $n\geq i-1$
the type $tp(\bar d_i/B_n)$ is principal,
and, therefore $tp(d_i/B_n)$ is also principal.
\end{minipage}
\end{equation}
Note that, by construction, $d_i$ has a principal type over every superset of $B_n$.

1) Let $\bar a=(a_1,a_2,...a_n)\in A$.
Since $a_j\in B\cup\underset{i<\omega}{\bigcup}D_i$ ($1\leq j\leq n$),
$\bar a$ is an enumeration of a set
$\{b_{i_1}, b_{i_2}, ..., b_{i_m}, d_{j_1}, d_{j_2}, ..., d_{j_l}\}$, where $m+l=n$.
Let $i=\max\{{i_1}, {i_2}, ..., {i_m}\}$,
$j=\max\{{j_1}, {j_2}, ..., {j_l}\}$, and
$k=\max\{i,j\}$.
By \eqref{eqS} we have that $tp(\bar d_j/B_k)$ is principal.
Then $tp(\bar b_i\bar d_j/B_k)$ is obviously principal,
what implies that $tp(\bar a/B_k)$ is principal as well.

2) Let $\mathfrak C\models T$ with $B\subseteq C$, and
let $p(\bar x)\in S_n(T)$ ($n<\omega$) be such that $p\in \mathcal D(\mathfrak A)$.
Since $p$ is realized in $\mathfrak A$ take an arbitrary realization
$\bar\alpha\in p(\mathfrak A)$.
By 1), there is $\bar b\in A$ such that
$tp(\bar\alpha/\bar b)$ is principal.
This type has to be realized in $\mathfrak C$.
Since $p=tp(\bar\alpha/\bar b)\restriction\emptyset$, it is also realized in the model
$\mathfrak C$ and is in its finite diagram.
And so $\mathcal D(\mathfrak A)\subseteq\mathcal D(\mathfrak C)$.

3) follows form 2).

\hfill $\Box$

Later in the article models obtained by the construction above are denoted as $\mathfrak A^B$.
%We use the notation $\mathfrak A^{\mathfrak M}$ to indicate construction over
%the universe $M$ of a structure $\mathfrak M$.

For natural $i\geq 1$ let $D^i$ be the set of elements of $\mathfrak A^B$
that were chosen during construction in
Theorem~\ref{Construction} as realizations of formulas from $\underset{j\leq i}{\bigcup}\Phi^j$.
When applied to the universe of a model $\mathfrak M$, the obtained structure $\mathfrak A^ M$
will coincide with $\mathfrak M$ itself, and each $A_n:=B_n\cup D^n$ will be a universe of a prime model
$\mathfrak A_n$ over the finite set $B_n$.
It is easy to see that $\mathfrak A_n\prec\mathfrak A_{n+1}\prec\mathfrak A^B$ for each $n<\omega$.
This way, Theorem~\ref{Construction} supports and proceeds S.V.~Sudoplatov's construction from \cite{Sud2018}
(Part 1, Theorem 1.1.3.1 in Russian version).
This construction can be be applied to a finite set $B$ as well,
and it will result in a prime model of $T$ over $B$.

\begin{corollary}\label{ConstructionCorollary1}
Let $\mathfrak M$ be a countable model of a small countable complete theory $T$,
and let $B_1$ and $B_2$ be subsets of $M$.

1) If $B_1\subseteq B_2\subseteq M$, then $\mathcal D(\mathfrak A^{B_1})\subseteq \mathcal D(\mathfrak A^{B_2})$;

2) If $\mathcal D(B_1)\subseteq \mathcal D(B_2)$, then
$\mathcal D(\mathfrak A^{B_1})\subseteq \mathcal D(\mathfrak A^{B_2})$.
\end{corollary}

\noindent{}{\it Proof of Corollary~\ref{ConstructionCorollary1}.}
1) Let $p\in \mathcal D(\mathfrak A^{B_1})$.
This means there exists $\bar b\in p(\mathfrak A^{B_1})$.
By Theorem~\ref{Construction}, 3), there exists $\bar b_1\in B_1$ such that $tp(\bar b/ \bar b_1)$ is principal.
Since $\bar b_1\in B_1\subseteq B_2\subseteq A^{B_2}$, $tp(\bar b/ \bar b_1)$ should be realized in $\mathfrak A^{B_2}$.
Therefore, $p=tp(\bar b/\bar b_1)\restriction \emptyset$ is realized in $\mathfrak A^{B_2}$ as well,
and is in its finite diagram, $\mathcal D(\mathfrak A^{B_2})$.

2) Let $p\in \mathcal D(\mathfrak A^{B_1})$. This means there exists $\bar b\in p(\mathfrak A^{B_1})$.
By Theorem~\ref{Construction}, 3), there exists $\bar b_1\in B_1$ such that $tp(\bar b/ \bar b_1)$ is principal.
Denote $p_1(x,\bar b_1):=tp(\bar b/ \bar b_1)$, and let $\varphi(\bar x,\bar b_1)$ be the isolating formula of this type.
We have that for every $P(\bar x, \bar b_1)\in tp(\bar b/ \bar b_1)$,
$\forall x\big(\varphi(\bar x,\bar y)\rightarrow P(\bar x,\bar y)\big)\in tp(\bar b_1)$.

Since $\bar b_1\in B_1\subseteq M$, $tp(b_1)\in\mathcal D(B_1)\subseteq \mathcal D(B_2)$.
Therefore there exists a model $\mathfrak M_2\models T$ and an element $\bar b_2\in B_2$ from $\mathfrak M_2$,
such that $tp(\bar b_1)=tp(\bar b_2)$.
Let  $p_1(\bar x,\bar b_2):=\{ P(\bar x,\bar b_2)\ |\ P(\bar x,\bar b_1)\in p(\bar x,b_1)\}$.
Then for every $P(\bar x, \bar b_2)\in tp(\bar b/ \bar b_2)$,
$\forall x\big(\varphi(\bar x,\bar y)\rightarrow P(\bar x,\bar y)\big)\in tp(\bar b_1)=tp(\bar b_2)$.
Analogically, $\exists x P(\bar x,\bar b_2)\in tp(\bar b_2)$.
By the statements above, $p_1(\bar x,\bar b_2)$ is a principal type over $\bar b_2$.
Then $p_1(\bar x,\bar b_2)$ is realized in $\mathfrak A^{B_2}$, and
$p_1(\bar x,\bar b_2)\restriction\emptyset=p\in\mathcal D(\mathfrak A^{B_2})$.

\hfill $\Box$

The following corollary shows that the construction in Theorem~\ref{Construction} preserves non-ho\-mo\-gen\-eity.

\begin{corollary}\label{homogeneous1}
Let $T$, $\mathfrak M$, and $B$ be as in Theorem~\ref{Construction},
and let $\bar a_1, \bar a_2, c\in B$ be such that
$tp(\bar a_1)=tp(\bar a_2)$, and
$tp(c\bar a_1)\neq tp(d \bar a_2)$ for every $d\in M$.
Then the structure $\mathfrak A=\mathfrak A^B$ is non-homogeneous.
\end{corollary}

\noindent{}{\it Proof of Corollary~\ref{homogeneous1}.}
To obtain a contradiction, suppose that there exists $d\in A$ such that  $tp(c\bar a_1)=tp(d\bar a_2)$.
Then $d$ should be either in $B$, which would contradict to the statement of the theorem, or
in $D_i$ for some $i<\omega$.
In the last case, we have that $tp(d/\bar b_n)$ is principal for all $n\geq k$ for a sufficient $k$.
But then we can find $m\geq n$ such that $\bar a_2\subseteq \bar b_m$ and the type $tp(d/\bar b_m)$ is isolated.
Therefore there is $d'\in M$ satisfying the type $tp(d/\bar b_m)\supseteq tp(d/\bar a_2)$ and such that
$tp(d'\bar a_2)=tp(d\bar a_2)=tp(c\bar a_1)$, what is a contradiction.

\hfill $\Box$

\begin{definition}\cite{Shelah1978}
Let $\mathfrak N$ be an $|A|^+$-saturated structure, where $A\subseteq N$, and let $p(\bar x),q(\bar y)\in S(A)$.
The type $p$ is {\bf weakly orthogonal} to the type $q$, $p\perp^wq$, if
$p(\bar x)\cup q(\bar y)$ is a complete $(ln(\bar x)+ln(\bar y))$-type.
\end{definition}

The following is a consequence of Corollary~\ref{homogeneous1}.
It is a part of Theorem~4 from \cite{BTYZ2015}.

\begin{corollary}\label{homogeneous2}
Let $\mathfrak M$ be a countable non-homogeneous model of be a small countable complete theory $T$,
and $p\in S(T)$ be such that $p\perp^w q$ for every $q\in \mathcal D(\mathfrak M)$.
Let $\bar a$ be a realization of $p$ in some elementary extension of $\mathfrak M$.
Then $\mathfrak A^{M\cup\{\bar a\}}$ is non-homogeneous.
\end{corollary}

\section{Quasi-successor formulas}\label{Squasi}

\begin{definition}\cite{B, BK2006, ABZ2014}
Let $\mathfrak M$ be a linearly ordered structure,
$A\subseteq M$, $\mathfrak M$ be $|A|^+$-saturated, and $p\in S_1(A)$ be non-algebraic.

1) An $A$-definable formula $\varphi(x, y)$ is said to be {\bf p--preserving (p-stable)} if there exist
 $\alpha$, $\gamma_1$, $\gamma_2 \in p(\mathfrak M)$ such that
$p(\mathfrak M)\cap \big(\varphi(M, \alpha)\setminus
\{\alpha\}\big)\ne\emptyset$ and $\gamma_1 < \varphi(M, \alpha) < \gamma_2$.

2) A $p$-preserving formula $\varphi(x, y)$ is said to be {\bf convex to the right (left)}
on $p$ if there exists
$\alpha$ $\in$ $p(\mathfrak M)$ such that $p(\mathfrak M)\cap \varphi(M, \alpha)$ is convex,
$\alpha$ is the left (right) endpoint of the set $\varphi(M, \alpha)$, and $\alpha\in \varphi(M, \alpha)$.

3) A $p$-preserving convex to the right (left) formula $\varphi(x,y)$ is a {\bf quasi-successor} on $p$
if for every $\alpha\in p(\mathfrak M)$ and every $\beta\in \varphi(M, \alpha)\cap p(\mathfrak M)$
\begin{center}
$p(\mathfrak M)\cap[\varphi(M, \beta)\backslash \varphi(M, \alpha)]\ne\emptyset$.
\end{center}
\end{definition}

For a quasi-successor formula $\varphi(x,y)$ we denote

\begin{tabular}{lcl}
$\varphi^0(x,y)$&$:=$&$(x=y)$;\\
$\varphi^1(x,y)$&$:=$&$\varphi(x,y)$;\\
$\varphi^{-1}(x,y)$&$:=$&$\varphi(y,x);$\\
$\varphi^n(x,y)$&$:=$&$\exists x_1, x_2,...,x_{n-1}
\big(\varphi(x_1,y)\land\varphi(x_2,x_1)\land...\land\varphi(x,x_{n-1})\big)$;\\
$\varphi^{-n}(x,y)$&$:=$&$\exists y_1, y_2,...,y_{n-1}
\big(\varphi(y,y_1)\land\varphi(y_1,y_2)\land...\land\varphi(y_{n-1},x)\big)$.
\end{tabular}

Let $\varphi(x,y)$ be a quasi-successor on a type $p(x)\in S_1(A)$.
For $\alpha,\beta\in p(\mathfrak M)$ we denote
\begin{center}
\begin{tabular}{ccl}
$V_{p(\mathfrak M),\varphi}(\alpha)$&$:=$&$\{\gamma\in p(\mathfrak M)\ |\ \exists n\in\mathbb Z, \gamma\in
\varphi^n(M, \alpha)\cap p(\mathfrak M)\}$;\\

$(V_{p,\varphi}(\alpha),V_{p,\varphi}(\beta))_{p(\mathfrak M)}$&$:=$&$\{\gamma\in p(\mathfrak M) \ |\
V_{p(\mathfrak M),\varphi}(\alpha)<\gamma<V_{p(\mathfrak M),\varphi}(\beta)\}$.
\end{tabular}
\end{center}

 When no confusion appears, we omit $\mathfrak M$ from the index: $V_{p,\varphi}(\alpha)$,
 $(V_{p,\varphi}(\alpha),V_{p,\varphi}(\beta))_{p}$.
We call the set $V_{p,\varphi}(\alpha)$ a {\bf neighborhood} (voisinage in French)
 of $\alpha$ in the 1-type $p$.

\begin{proposition}\label{Ext}
Let T be a small countable complete theory, $p\in S_n(T)$.
If $I(T\cup p(\bar c),\aleph_0)=2^{\aleph_0}$, where $\bar c$ is a tuple of new constants, then $I(T,\aleph_0)=2^\aleph_0$.
\end{proposition}

The theory $T'=T\cup p(\bar c)$ is called an {\bf inessential expansion} of the theory $T$.
The Proposition~\ref{Ext} allows us to transfer to an inessential expansion of a theory
in the case of maximality of the number of countable non-isomorphic models.
It holds because every model of $T'$ is a model of $T$ when the tuple $\bar c$ is forgotten,
and every countable model of $T$ can produce at most $\aleph_0$ countable non-isomorphic models
of $T'$ since the set of all realizations of $p$ is no more than countable.

\begin{theorem}\cite{ABZ2014}\label{abzQuasi}
Let $T$ be a linearly ordered theory,
$A$ be a finite subset of a countable saturated model $\mathfrak M\models T$,
and let $p(x)\in S_1(A)$ be non-algebraic.
If there exists an  $A$-definable quasi-successor on $p$, then $T$ has $2^{\aleph_0}$ countable models.
\end{theorem}

\noindent{}{\it Proof of Theorem~\ref{abzQuasi}.}
Since every non-small theory has $2^{\aleph_0}$ countable non-isomorphic models,
it remains to prove the case, when the theory $T$ is small.

For simplicity we extend our language to $\mathcal L(A)$ and work in the theory
$T\cup tp(\bar a)$, where $\bar a$ is an enumeration of the set $A$.
Without loss of generality assume that the quasi-successor formula $\varphi$
is convex to the right on $p$.

We will prove the existence of an infinite indiscernible sequence of elements
from the set of realizations of $p$ such that
every two realizations located between neighborhoods of its elements
have the same convex 1-type over every its finite subsequence.

Consider the following (not necessary compete) type
$q_0(x,y):=\big\{x<y\big\}\cup p(x)\cup p(y)\cup
\big\{\neg\varphi^n(y,x)\ |\ n\in\mathbb N\big\}\cup
\Big\{R(y,x)\ |\ R(y,x)$ is a convex to the right on $p$ formula such that
for all $n\in\mathbb N$ and all $\alpha\in p(\mathfrak M)$,
$\big((\varphi^n(M, \alpha)\cap p(\mathfrak M))\subseteq R(M, \alpha)\big)\Big\}\cup
\Big\{L(x,y)\ |\ L(x,y)$ is a convex to the left on $p$ formula such that
for all $n\in\mathbb N$ and all $\beta\in p(\mathfrak M)$,
$\big((\varphi^{-n}(M, \beta)\cap p(\mathfrak M))\subseteq L(M, \beta)\big)\Big\}$.
Consistence of $q_0(x,y)$ is verified directly.
Notice that, when $T$ is weakly o-minimal, $q_0$ is a complete 2-type.

\begin{lemma}\label{ConvexTypesQuasi}
There exists a complete 2-type $q(x,y)\supseteq q_0(x,y)$ such that for every tuple $(\alpha, \beta)$
realizing $q(x,y)$, for all
$\gamma_1,\gamma_2\in (V_{p,\varphi}(\alpha),V_{p,\varphi}(\beta))_{p(\mathfrak M)}$,
\begin{center}
$tp^c(\gamma_1/\alpha\beta)=tp^c(\gamma_2/\alpha\beta)$.
\end{center}
\end{lemma}

\noindent{}{\it Proof of Lemma~\ref{ConvexTypesQuasi}.}
We  call the  conclusion of Lemma~\ref{ConvexTypesQuasi} a {\bf $V_{p,\varphi}$-$2$-condition}.
Towards a contradiction suppose that the conclusion of Lemma~\ref{ConvexTypesQuasi}
is not true. This assumption leads to the existence of a continuum of 1-types.
Namely, we prove existence of a formula separating two neighborhoods of two elements of the type,
that simultaneously separates the neighborhood of some element of the type.

Thus for every 2-type $q(x,y)\supseteq q_0(x,y)$ and for $(\alpha, \beta)$ realizing $q$ in $\mathfrak M$  there are
$\gamma_1,\gamma_2\in (V_{p,\varphi}(\alpha),V_{p,\varphi}(\beta))_{p(\mathfrak M)}$
and a convex formula $H'(z,\alpha,\beta)$ such that
$\gamma_1\in H'(M,\alpha,\beta)<\gamma_2$.
For convenience, replace $H'$ with the formula determining the set
$(H'(M,\alpha,\beta)^+)^-$.
Denote $H(z,x,y):=\exists z_1\exists z_1(x<z_1<z_2<y\land H'(z_1,x,y)\land \neg H'(z_2,x,y))$.

Notice that for any 2-formula $Q(x,y)\in q$, $Q(\alpha,M)\cap V_{p,\varphi}(\alpha)\not =\emptyset$.
If we suppose that this intersection is empty,
then the $\alpha $-formula $y\geq\alpha \land Q(\alpha, y)^-$ is convex to the right and has to belong to $q(\alpha,y)$.
A contradiction.

Denote $C(\alpha, y):=\exists z H(z,\alpha,y)$.
Since $\mathfrak M\models C( \alpha, \beta)$,
and, consequently, $C(\alpha, y)\in q(\alpha, y)$,
 there is $b \in C( \alpha, M)\cap V_{p,\varphi}(\alpha). $
Let $n<\omega$ be such that  $H(M,\alpha,b)\cap \varphi^n(M,\alpha)\not = \emptyset $  and
$\neg H(M,\alpha,b)\cap \varphi^n(M,\alpha)\not = \emptyset  $.
Then there are two elements $c_1,c_2$
from $ p(\mathfrak M)$ such that $\mathfrak M\models \big(H(c_1,\alpha,b)\land\neg H(c_2,\alpha,b)\land\varphi(c_2,c_1)\big)$.
Thus, the 2-formula
\[S (x,y):= \exists z_1\exists z_2 \big(H(z_1,x,y)\land\neg H(z_2,x,y)\land\varphi(z_2,z_1)\big)\]
belongs to $q(x,y)$ because any element from $V_{p,\varphi}(\alpha)$ that satisfies $C(\alpha, y)$
satisfies $S(\alpha, y)$.

For  $Q(x,y)\in q$,
$P(x)\in p$, $n,k,l< \omega$ such that $k+l< n$ introduce the following
2-formula:

\begin{center}
$\Gamma[Q,P,n,k,l](x,y):=\big(x<y\land Q(x,y)\land S(x,y)\land\neg\varphi^n(y,x)\big)\rightarrow
\exists z_1\exists z_2\big(x<z_1<z_2<y\land P(z_1)\land P(z_2)\land \varphi(z_2,z_1)\land\neg\varphi^k(z_1,x)
\land\neg\varphi^l(y,z_2)\land H(z_1,x,y)\land\neg H(z_2,x,y)\big)$.
\end{center}

Let $\langle Q_j \rangle_{j<\omega}$ and $\langle P_i\rangle_{i<\omega}$ be two strictly decreasing  sequences
of 2-formulas and 1-formulas such that  for every $j, i<\omega$
\begin{center}
$\mathfrak M \models \forall x \forall y \big(Q_{j+1}(x,y)\to Q_j(x,y)\big)\land \forall x\big(P_{i+1}(x)\to P_i(x)\big)$
\end{center}
and
\begin{center}
$q(\mathfrak M)= \underset{j<\omega}\bigcap  Q_j (M^2)$, $p(\mathfrak M)= \underset{i<\omega}\bigcap P_i(M)$.
\end{center}

\begin{claim}\label{CS}
There are five increasing sequences of natural numbers:
$\langle j(m), i(m), n(m), k(m),$ $l(m)\rangle_{m<\omega}$
such that for every $m<\omega$
\begin{center}
$\mathfrak M \models \forall x \forall y \Gamma[Q_{j(m)},P_{i(m)},n(m),k(m),l(m)](x,y)$.
\end{center}
\end{claim}

{\it Proof of Claim~\ref{CS}.}
Firstly, towards a contradiction, suppose that for every
$Q(x,y)\in q, P(x)\in p, n,k,l<\omega$, $n> k+l$,
$\mathfrak M\models \exists x\exists y\neg\Gamma[Q,P,n,k,l](x,y)$, that is,
\begin{center}
$\mathfrak M \models\exists x\exists y \Big[x<y\land Q(x,y)\land S(x,y)\land \neg\varphi^n(y,x)\land
\forall z_1\forall z_2\Big(\big(x<z_1<z_2<y\land P(z_1)\land P(z_2)\land \varphi(z_2,z_1)\land\neg\varphi^k(z_1,x)
\land\neg\varphi^l(y,z_2)\land H(z_1,x,y)\big)\rightarrow H(z_2,x,y)\Big)\Big]$
\end{center}
 or, equivalently,
\begin{center}
$\mathfrak M \models\exists x\exists y \Big[x<y\land Q(x,y)\land S(x,y)\land \neg\varphi^n(y,x)\land
\forall z_1\forall z_2\Big(\big(x<z_1<z_2<y\land P(z_1)\land P(z_2)\land \varphi(z_2,z_1)\land\neg\varphi^k(z_1,x)\land\neg\varphi^l(y,z_2)\land
\neg H(z_2,x,y)\big)\rightarrow \neg H(z_1,x,y)\Big)\Big]$.
\end{center}

Denote
\begin{center}
$\Psi[Q,P,n,k,l](x,y):= \neg\Gamma[Q,P,n,k,l](x,y)\land \exists z\big(P(z)\land \neg\varphi^k(z,x)\land x<z<y \land H(z,x,y)\big)$;

$\Theta[Q,P,n,k,l](x,y):= \neg\Gamma[Q,P,n,k,l](x,y)\land \exists z\big(P(z)\land \neg\varphi^l(y,z)\land x<z<y \land \neg H(z,x,y)\big)$.
\end{center}
From definitions of these formulas we have
\begin{equation}
\label{eq1}
\mathfrak M\models \forall x\forall y\big[\neg\Gamma[Q,P,n,k,l](x,y)\leftrightarrow \big(\Psi[Q,P,n,k,l](x,y)\lor\Theta[Q,P,n,k,l](x,y)\big)\big].
\end{equation}

{\bf Step 0.} Fix arbitrary $l_0<\omega $, $i_0<\omega$.

Consider the set of 2-formulas
\begin{center}
$r(i_0,l_0)(x,y):=\{\Psi[Q_j,P_{i_0},n,k,l_0](x,y)\ |\ j<\omega, \; n< \omega, \; k+l_0<n \}$.
\end{center}
Notice that $\mathfrak M\models \forall x\forall y\big(\Psi[Q_{j+1},P_{i_0},n+1,k+1,l_0](x,y) \to \Psi[Q_j,P_{i_0},n,k,l_0](x,y)\big)$.
If the set $r(i_0,l_0)(x,y)$ is consistent, then its realization $(\alpha',\beta')$ has the property:
\\ $H(M, \alpha',\beta')\cap V_{p,\varphi}(\beta')\not =\emptyset.$
Contradiction with definition of the 2-type $q(x,y)$. Thus, for some $j_0 <\omega$, $n_0<\omega$, $k_0< \omega$
\begin{equation}
\label{eq2}
\mathfrak M\models \forall x\forall y\big(\neg\Psi[Q_{j_0},P_{i_0},n_0,k_0,l_0](x,y)\big).
\end{equation}

Then by \eqref{eq1}
\begin{center}
$\mathfrak M\models \exists x\exists y\big(\Theta[Q_{j_0},P_{i_0},n_0,k_0,l_0](x,y)\big)$.
\end{center}

Consider the set of 2-formulas  $r(i_0,k_0)(x,y):=$
\begin{center}
$\{\Theta[Q_j,P_{i_0},n,k_0,l](x,y)\
|\ j_0<j<\omega, \; n_0 <n< \omega, \; k_0+l <n, l_0 <l \}$.
\end{center}

If this set is consistent then we obtain that any realization of this set of 2-formulas gives the property:
$\neg H(M, \alpha',\beta')\cap V_{p,\varphi}(\alpha)\not =\emptyset$.
The last contradicts to the definition of $H$ for realization of $q(x,y)$.
Thus, there exist $j_1, n_1, l_1<\omega$ such that
\begin{equation}
\label{eq3}
\mathfrak M\models \forall x\forall y(\neg\Theta[Q_{j_1},P_{i_0},n_1,k_0,l_1](x,y)).
\end{equation}

Considering properties of these formulas and \eqref{eq2} and \eqref{eq3} we obtain
\begin{center}
$\mathfrak M\models \forall x\forall y\Gamma[Q_{j_1},P_{i_0},n_1,k_0,l_0](x,y)$.
\end{center}

Denote $j(0):=j_1$, $i(0):=i_0$, $n(0):=n_1$, $k(0):=k_0$, $l(0):=l_0$.
Then

\begin{center}
$\mathfrak M\models \forall x\forall y\Gamma[Q_{j(0)},P_{i(0)},n(0),k(0),l(0)](x,y)$.
\end{center}

{\bf Step m+1.} We repeat all considerations of Step 0 with complimentary conditions, namely:
$l_0> l(m)$, $i_0> i(m)$, $j_0> j(m)$, $k_0> k(m)$, $n_0> n(m)$.

\hfill$\Box$ Claim~\ref{CS}

The formula $H(x,\alpha,y)$ is increasing,
$\mathfrak M\models\forall x \forall y \Big[\exists z H(z,x,y) \to \Big(\forall y'\big((\exists z' H(z',x,y') \land y<y') \to
\forall z'(H(z',x,y)\to H(z', x,y'))\big)\Big)\land  \exists y''\exists z'' \Big(y<y''\land
H(z'',x,y'')\land \neg H(z'',x,y)\Big)\Big]$.
When $m$ approaches infinity, $n(m)$, $k(m)$, $l(m)$ and the
difference $(n(m)-l(m))$ also approach infinity.
It means that
\begin{equation}
\label{eq4}
\begin{split}
\mathfrak M \models &\forall x \forall y [\exists z_1 \exists z_2 (H(z_1,x,y)\land \varphi(z_2,z_1) \land\neg H(z_2,x,y))\to\\
&\forall y_1 ((\varphi(y_1,y)  \land \exists z_3 H(z_3,x,y_1 )) \to \forall z_4(z_4 > \varphi (M, z_1) \to \neg H(z_4,x, y_1)))].
\end{split}
\end{equation}

Thus, by Claim~\ref{CS}, we have three neighborhoods:
\begin{equation}
\label{eq5}
V_{p,\varphi}(\alpha)<V_{p,\varphi}(\delta)< V_{p,\varphi}(\beta),\
p(\mathfrak M)\cap H (M,\alpha,\beta)\cap V_{p,\varphi}(\delta)\not = \emptyset
\end{equation}
for some $\delta\in p(\mathfrak M)$, that is determined by $H(x,\alpha,\beta)$ and because $S(x,y)\in q$.

We will construct a countable number of convex disjoint 1-$\alpha\beta$-formulas
 by $K$-steps by using the family of 1-$\alpha\beta$-formulas of $G$-kind.

Namely, {\it for every $n<\omega$, for every $\tau \in 2^n$
(here, $2^n$ is the set of all bit strings of length $n$) we define convex
1-$\alpha\beta$-formulas $G_{\tau}$ and $K_{\tau}$
such that there are $\mu, \gamma \in\underset{k<n}{\bigcup}2^k$ with
\begin{equation}
\tag{$5_\tau$}
V_{p,\varphi}(\delta_{\mu})< V_{p,\varphi}(\delta_{\tau})< V_{p,\varphi}(\delta_{\gamma}),
\end{equation}
or, equivalently, since $G_{\tau}(M)\cap p(\mathfrak M) \subset V_{p,\varphi}(\delta_{\tau})$,
$K_{\tau}(x):= G_{\mu}(M)< x < G_{\gamma}(M)$
\begin{equation}
\tag{$5_\tau$}
G_{\mu}(M)< G_{\tau}(M)< G_{\gamma}(M),
\end{equation}
\begin{equation}
\tag{$6_\tau$}
K_{\tau0}(M)< G_{\tau}(M)< K_{\tau1}(M)
\end{equation}
and for any initial substring $\mu$ of $\tau$
($\mu \subset \tau$),
\begin{equation}
\tag{$7_\tau$}
K_{\tau}(M)\subset K_{\mu}(M),
\end{equation} or, equivalently,
\begin{equation}
\tag{$7_\tau$}
K_{\mu0}(M)\cup K_{\mu1}(M)\subset K_{\mu}(M).
\end{equation}}
Here and below, $\mu0$, $\mu1$, ... stand for concatenations, $\mu^\frown0$, $\mu^\frown1$, ... .
Notice that for $\mu,\gamma$ we suppose that either $\mu0=\tau$ or $\gamma1=\tau$:
if $\mu0=\tau$, then the length of $\gamma$ is less than $n-1$, and
if $\gamma1=\tau$, then the length of $\mu$ is less than $n-1$.

There are 2 cases:

$a)$ $tp (\delta/\alpha)=q(\alpha, y)$ and $tp(\delta/\beta )=q(x,\beta)$;

$b)$  $tp (\delta/\alpha)\not =q(\alpha, y)$ or $tp(\delta/\beta )\not =q(x,\beta)$.

$a)$ In this case, for any  $m< \omega$ there are $n<\omega $, $m_1<\omega$  ($m_1>m$), such that
\begin{equation}
\tag{8}
\begin{split}
\mathfrak M \models & \forall x \forall y \big[\big(Q_{m_1}(x,y)\land \neg \varphi^n(y,x)\big) \to \\
& \exists z \exists z_1 \big(Q_m(x,z) \land Q_m(z,y) \land \varphi(z_1,z)
\land H(z,x,y)\land \neg H(z_1,x,y)\big)\big].
\end{split}
\end{equation}

{\bf Step K-$0$.}

Denote $G(x,\alpha, \beta):=
H(x,\alpha, \beta)\land\exists y\big(\varphi(y,x)\land \neg H(y,\alpha,\beta)\big)$.
As we already showed by Claim~\ref{CS}, $G(M, \alpha, \beta)\cap p(\mathfrak M)\neq\emptyset$, and
by \eqref{eq4} for some
$\delta\in(V_{p,\varphi}(\alpha),V_{p,\varphi}(\beta))$, $G(M, \alpha, \beta)\cap p(\mathfrak M)\subseteq V_{p,\varphi}(\delta)$.
Thus,
\begin{equation}
\tag{$5_0$}
V_{p,\varphi}(\alpha)< V_{p,\varphi}(\delta) <V_{p,\varphi}(\beta).
\end{equation}
This means we have three disjoint convex 1-$\alpha\beta$-subformulas
of $K_{\emptyset}:= \alpha <x < \beta$: $G(M,\alpha,\beta)$,
 $K_0 (x,\alpha, \beta):= \alpha <x < G(M,\alpha,\beta)$
and $K_1(x,\alpha,\beta):=  G(M,\alpha,\beta)<x<\beta$.
Thus,
\begin{equation}
\tag{$6_0$}
K_0(M)<G(M)<K_1(M),
\end{equation}
\begin{equation}
\tag{$7_0$}
K_0(M)\cup K_1(M)\subset K_{\emptyset}(M).
\end{equation}
{\bf Step K-$1$.}
Let
\begin{center}
$H_0(x,\alpha,\beta):=\exists z(H(x,\alpha,z)\land G(z,\alpha,\beta))$;
$H_1(x,\alpha,\beta):=\exists z(H(x,z,\beta)\land G(z,\alpha,\beta))$.
\end{center}

These 1-$\alpha \beta$-formulas are defined correctly because $tp(\alpha \delta)=tp (\delta \beta)=q(x,y)$.

Then by~\eqref{eq4}
\begin{equation}
\tag{$8_1$}
V_{p,\varphi}(\alpha)\subset H_0(M,\alpha,\beta)<V_{p,\varphi}(\delta)\subset H_1(M,\alpha,\beta)<V_{p,\varphi}(\beta).
\end{equation}

Now denote

$G_0(x,\alpha,\beta):= H_0(x,\alpha,\beta)\land \exists y (\neg H_0(y,\alpha,\beta)\land \varphi(y,x))$.

$G_1(x,\alpha,\beta):= H_1(x,\alpha,\beta)\land \exists y (\neg H_1(y,\alpha,\beta)\land \varphi(y,x))$.

By \eqref{eq4} and \eqref{eq5} there are
$\delta_0\in(V_{p,\varphi}(\alpha),V_{p,\varphi}(\delta))_p$ and
$\delta_1\in(V_{p,\varphi}(\delta),V_{p,\varphi}(\beta))_p$ such that
\begin{equation}
\tag{9}
tp(\alpha,\delta_0)=tp(\delta_0 , \delta )= tp(\delta,\delta_1)= tp(\delta_1,\beta)=q(x,y)
\end{equation}
and $G_0(M,\alpha,\beta)\cap p(\mathfrak M)\subset V_{p,\varphi}(\delta_0)$;
$G_1(M,\alpha,\beta)\cap p(\mathfrak M)\subset V_{p,\varphi}(\delta_1)$.
Thus,
\begin{center}
$V_{p,\varphi}(\alpha)<V_{p,\varphi}(\delta_0)<V_{p,\varphi}(\delta)<V_{p,\varphi}(\delta_1)<V_{p,\varphi}(\beta)$.
($5_0$) ($5_1$)
\end{center}
This means
\begin{center}
$ V_{p,\varphi}(\alpha)< G_0(M,\alpha,\beta)< G(M, \alpha,\beta)< G_1(M,\alpha,\beta)<V_{p,\varphi}(\beta)$.
\end{center}

Let $K_{00}(x,\alpha,\beta):= \alpha < x < G_0(M,\alpha,\beta)$, \\
$K_{01}(x,\alpha,\beta):= G_0(M,\alpha,\beta)< x <  G(M, \alpha,\beta)$,\\
$K_{10}(x,\alpha,\beta):= G(M, \alpha,\beta)< x< G_1(M,\alpha,\beta)$, and\\
$K_{11}(x,\alpha,\beta):=G_1(N,\alpha,\beta) < x < \beta$.
Then

\begin{center}
$\alpha <  K_{00}(M)< K_{01}(M)< G(M,\alpha,\beta)< K_{10}(M)< K_{11}(M)<\beta$,
\end{center}
\begin{equation}
\tag{$6_0$}
K_{00}(M)<G_0(M)< K_{01}(M),
\end{equation}
\begin{equation}
\tag{$6_1$}
K_{10}(M)< G_1(M)<K_{11}(M),
\end{equation}
\begin{equation}
\tag{$7_0$}
K_{00}(M, \alpha,\beta)\cup K_{01}(M, \alpha,\beta)\subset K_{0}(M, \alpha,\beta),
\end{equation}
\begin{equation}
\tag{$7_1$}
K_{10}(M, \alpha,\beta)\cup K_{11}(M, \alpha,\beta)\subset K_{1}(M, \alpha,\beta).
\end{equation}

{\bf Step K-$(n+1)$.} By induction supposition for $\tau \in 2^n$
there are $\delta_{\tau}\in p(\mathfrak M)$ and a 1-$\alpha\beta$-formula
$G_{\tau}$ such that $ (p(\mathfrak M)\cap G_{\tau}(M))\subset V_{p,\varphi}(\delta_{\tau})$
and there are $\mu,\gamma\in \cup_{k<n}2^k$ such that
\begin{center}
$V_{p,\varphi}(\delta_{\mu})< V_{p,\varphi}(\delta_{\tau})< V_{p,\varphi}(\delta_{\gamma})$,
\end{center}
and, consequently, $G_{\mu}(M)< G_{\tau}(M)< G_{\gamma}(M).$

By using the $3$-formula $H$ that provides Claim~\ref{CS} we introduce the convex 1-$\alpha\beta$-formulas:

\begin{equation}
\tag{$10_{\tau0}$}
H_{\tau0}(x):=\exists z_1\exists z_2 \big(G_{\mu}(z_1)\land G_{\tau}(z_2)\land H(x,z_1,z_2)\big);
\end{equation}

\begin{equation}
\tag{$10_{\tau1}$}
H_{\tau1}(x):=\exists z_1\exists z_2 \big(G_{\tau}(z_1)\land G_{\gamma}(z_2)\land H(x,z_1,z_2)\big).
\end{equation}

Notice that by assuming that ALL complete 2-types containing $q_0(x,y)$ do
not satisfy the condition of Lemma~\ref{ConvexTypesQuasi} have to have such $3$-formula $H$  provided by Claim~\ref{CS}.
For the case b) we denote this $3$-formula by $H^{\mu\tau}(z,x,y)$,
and denote $q^{\mu\tau}(x,y):=tp(\delta_{\mu}, \delta_{\tau})$.

\begin{center}
$V_{p,\varphi}(\delta_{\mu})\subset H_{\tau0}(M) < V_{p,\varphi}(\delta_{\tau})\subset H_{\tau1}(M)
< V_{p,\varphi}(\delta_{\gamma})$,
\end{center}

\begin{center}
$V_{p,\varphi}(\delta_{\mu})<V_{p,\varphi}(\delta_{\tau0}) < V_{p,\varphi}(\delta_{\tau})<V_{p,\varphi}(\delta_{\tau1})< V_{p,\varphi}(\delta_{\gamma})$. ($5_{\tau0}$), ($5_{\tau1}$)
\end{center}

Define $G_{\tau0}(x)$ and $G_{\tau1}(x)$:

$G_{\tau0}(x):= \exists z(H_{\tau0}(x)\land \varphi(z,x)\land \neg H_{\tau0}(z))$,
by $(5_{\tau0})$ $G_{\tau0}(M) \subset V_{p,\varphi}(\delta_{\tau0})$;

$G_{\tau1}(x):= \exists z(H_{\tau1}(x)\land \varphi(z,x)\land \neg H_{\tau1}(z))$
by $(5_{\tau1})$ $G_{\tau1}(M) \subset V_{p,\varphi}(\delta_{\tau1})$.

\begin{center}
$G_{\mu}(M)<G_{\tau0}(M) < G_{\tau}(M) <G_{\tau1}(M) < G_{\gamma}(M)$,
\end{center}
$K_{\tau0}(x):=G_{\mu}(M)<G_{\tau0}(M)< x < G_{\tau}(M)$ and $K_{\tau1}(x):=G_{\tau}(M) < x < G_{\gamma}(M)$. Then
\begin{equation}
\tag{$6_{\tau}$}
K_{\tau0}(M) < G_{\tau}(M)< K_{\tau1}(M).
\end{equation}
\begin{equation}
\tag{$7_{\tau}$}
K_{\tau0}(M)\cup K_{\tau1}(x)\subset K_{\tau}(M).
\end{equation}

Case $b)$ is considered like case $a)$, in definitions of $H_{\tau0}$ and $H_{\tau1}$ in ($10_{\tau0}$) and $(10_{\tau1})$
$H^{\mu\tau}$ and $H^{\tau\gamma}$ should be used.

Thus, for an arbitrary $\nu \in 2^{\omega}$ there is a consistent set of 1-$\alpha\beta$-formulas
$r_{\nu}:= \{K_{\nu(n)}(x)\ |\ n<\omega\}$.
This contradicts with $T$ being small.

\hfill $\Box$ Lemma~\ref{ConvexTypesQuasi}

A type $r\in S_1(A)$ is {\bf irrational} if the sets $r^c(\mathfrak M)^+$ and $r^c(\mathfrak M)^-$ are both non-definable.
By Lemma~\ref{ConvexTypesQuasi} there exists a 2-type $q(x,y)$ such that
for any $\alpha,\beta\in p(\mathfrak M)$ the
following holds:  if $tp(\alpha, \beta)=q(x,y)$, for any $\gamma \in (V_{p,\varphi}(\alpha),V_{p,\varphi}(\beta))_{p}$,
the type
$tp^c(\gamma/ \alpha\beta)(\mathfrak M)= (V_{p,\varphi}(\alpha)^+ < x < V_{p,\varphi}(\beta)^-) (\mathfrak M)$
is irrational and, therefore, $tp(\gamma/\alpha\beta)$ is non-principal.

It follows from $V_{p,\varphi}$-2-condition for $q(x,y)$, for any 3-formula $H(z,x,y)$ with the condition
$\mathfrak M \models \forall x \forall y \big(H(M,x,y) < y\land \forall z(H(z,x,y)\leftrightarrow H(z,x,y)^{+-})\big)$,
there are two cases:

1) There is  $k_0<\omega$ such that
$\forall z \big(\neg \varphi^{k_0}(z,x)\to \neg H (z,x,y)\big)\in q(x,y)$, or, equivalently, for
$T_1[H](x,y):=\forall x_1\Big(\varphi(x,x_1) \to  \exists z\big(H(z,x,y)\land \neg H(z,x_1,y)\big)\Big)
\land \forall y_1\Big(\varphi (y_1,y)\to \forall z\big(\neg H(z,x,y) \leftrightarrow \neg H(z,x,y_1)\big)\Big)$, 
$T_1 [H](x,y)\in q$.

2) There is  $l_0< \omega$ such that
$\forall z \Big(\big(x< z \land \neg \varphi^{l_0}(y,z)\big)\to H(z,x,y)\Big)\in q(x,y)$, or, equivalently, for
$T_2 [H](x,y):=\forall y_1\Big(\varphi(y_1,y) \to\exists z\big(H(z,x,y)\land\neg H(z,x,y_1)\big)\Big) \land
\forall x_1\Big(\varphi (x,x_1)\to \forall z\big(H(z,x,y) \leftrightarrow H(z,x_1,y)\big)\Big)$, $T_2 [H](x,y)\in q$.

Since $\forall x \forall y (T_1[H](x,y) \to \neg T_2 [H](x,y)) \in T$ and
$\forall x \forall y\big(T_2[H](x,y) \to \neg T_1 [H](x,y)\big)\in T$,
$\big(T_1[H](x,y) \lor T_2 [H](x,y)\big)\in q$.

Thus, the following set of 2-formulas $q_1(x,y):= q_0(x,y)\cup \Big\{T_1[H](x,y) \lor T_2 [H](x,y)\ |$
$\mathfrak M \models \forall x \forall y$$\Big(H(M,x,y) < y\land \forall z\big(H(z,x,y)
\leftrightarrow H(z,x,y)^{+-}\big)\Big) \Big\}$
is consistent, and every complete extension of $q_1$ satisfies $V_{p,\varphi}$-2-condition.
In fact, the proof of Lemma~\ref{ConvexTypesQuasi} is consistency of $q_1(x,y)$.

As a corollary of the proof of Lemma~\ref{ConvexTypesQuasi} we obtain the following lemma.

\begin{lemma}\label{MoreConvexTypesQuasi}
For every $n$, $1<n<\omega$, there is an $n$-type $q^n(x_1,\dots, x_n)$ that satisfies
{\bf $V_{p,\varphi}$-$n$-condition}:
for every increasing sequence
$\alpha_1,\alpha_2, \dots, \alpha_n \in p(\mathfrak M)$ with $tp(\alpha_1,\alpha_2, \dots, \alpha_n)=q$,
for every $i$ ($1\leq i<n$), and every $\gamma_1,\gamma_2 \in (V_{p,\varphi}(\alpha_i),V_{p,\varphi}(\alpha_{i+1}))_p$, $tp^c(\gamma_1/\bar \alpha)= tp^c(\gamma_2/\bar \alpha)$, and
\begin{center}
$tp^c(\gamma_1/\bar \alpha)(\mathfrak M)=(V_{p,\varphi}(\alpha_i)^+< x < V_{p,\varphi}(\alpha_{i+1})^-)_
{p(\mathfrak M)}(\mathfrak M)$,
\end{center}
where $\bar\alpha:=(\alpha_1,\alpha_2, \dots, \alpha_n)$.
\end{lemma}

\noindent{}{\it Proof of Lemma~\ref{MoreConvexTypesQuasi}}
Let $\bar x_n=\langle x_1,x_2,\dots, x_n\rangle$,
denote $x_n^i:= \langle x_1, x_2,\dots, x_{i-1}, x_{i+1}, \dots, x_n\rangle$.

By induction on $n$ ($1<n<\omega$) will prove the existence of $q^n(\bar x_n)$ with  $V_{p,\varphi}$-$n$-condition.
For $n = 2$ this is Lemma~\ref{ConvexTypesQuasi},
we assume existence of $q^{n-1}$ that satisfies $V_{p,\varphi}$-$(n-1)$-condition.

Let $q_0^n(x_1,\dots, x_n):=\underset{1\leq i< n}\bigcup\big(V_{p,\varphi}(x_i) < V_{p,\varphi}(x_{i+1})\big)
\cup \big\{R(x_{i+1},\bar x^{i+1}_n)\ |$ for all $k<\omega$
$\mathfrak M\models\forall z\big(\varphi^k(z,x_i)\to R(z,\bar x^{i+1}_n)\big)\big\}
\cup \big\{L(x_i,\bar x_n^{i})\ |$  for all $l<\omega,\ \mathfrak M\models
\forall z\big(\varphi^l(x_{i+1},z)\rightarrow L(z,\bar x_n^{i})\big)\big\}$.
From existence of $q^{n-1}$ with $V_{p,\varphi}$-$(n-1)$-condition follows consistence of $q_0$.

Towards a contradiction suppose that every complete extension of $q_0^n(\bar x)$ does not satisfy $V_{p,\varphi}$-$n$-condition.
Let $\bar \alpha$ be such that $q_0^n (\bar x) \subset tp(\bar \alpha)=q(\bar x)$.
Since $q$ does not satisfy $V_{p,\varphi}$-$n$-condition, for some $i$ ($1\leq i<n$)
there exists a convex 1-$\bar \alpha$-formula $H^{\bar \alpha}_i(z,\bar \alpha)$ such that
\begin{center}
$V_{p,\varphi}(\alpha_i)\subset H^{\bar \alpha}_i(M,\bar \alpha)<  V_{p,\varphi}(\alpha_{i+1})$.
\end{center}
Then, by using Claim~\ref{CS}, we obtain that  there are $\gamma, \gamma_1  \in p(\mathfrak M)$
such that $\gamma  \in H^{\bar \alpha}_i(M, \bar \alpha)$  and
$\mathfrak M \models \varphi (\gamma_1,\gamma)\land \neg  H^{\bar \alpha}_i(\gamma_1, \bar \alpha)$.
Denote by $G_{i,i+1}^{\gamma}(z,\bar \alpha):= \exists z_1(\varphi (z_1,z)\land
H^{\bar \alpha}_i(z, \bar \alpha) \land \neg  H^{\bar \alpha}_i(z_1,\bar \alpha))$.
Then $G_{i,i+1}^{\gamma}(M,\bar \alpha)\cap p(\mathfrak M)\subset V_{p,\varphi}(\gamma)$ and
$V_{p,\varphi}(\alpha_i)< V_{p,\varphi}(\gamma)< V_{p,\varphi}(\alpha_{i+1})$.
Thus, we have $(n + 1)$ neighborhoods instead of $n$ and the 1-$\bar \alpha$-formula such that
$\mathfrak M \models \forall z \big( G_{i,i+1}^{\gamma}(z,\bar \alpha)\to \exists z_1
(\varphi (z_1,z)\land \neg G_{i,i+1}^{\gamma}(z_1,\bar \alpha))\big)$.

\begin{claim}\label{CSS}
Let $q_0^n(\bar x)\subseteq tp(\bar \delta)$ and for any $k \;(1\leq k \leq n)$, there exists $G_k(z,\bar \alpha)$ such that
$\delta_k \in (G_k(M,\bar \alpha)\cap p(\mathfrak M)) \subset V_{p,\varphi}(\delta_k)$, and
$\mathfrak M \models \forall z \Big( G_k(z,\bar \alpha)\to \exists z_1\big(\varphi (z_1,z)\land \neg G_k(z_1,\bar \alpha)\big)\Big)$.

By our assumption $tp(\bar \delta)$ does not satisfy $V_{p,\varphi}$-$n$-condition.
Then for some $i$ ($1\leq i <n$) there exists a convex 1-$\bar \delta$-formula $H_i(z,\bar \delta)$ such that
\begin{center}
$V_{p,\varphi}(\delta_i)\subset H_i (M,\bar \delta)< V_{p,\varphi}(\delta_{i+1})$.
\end{center}
Then for $q'(x,y):=tp(\delta_i \delta_{i+1}/\bar \delta^{i,i+1})$, by using arguments for Claim~\ref{CS},
there exists  $\beta \in p(\mathfrak M)$ such that
$ H_i (M,\bar \delta)\cap V_{p,\varphi}(\beta)\not = \emptyset$ and
$\neg H_i (M,\bar \delta)\cap V_{p,\varphi}(\beta)\not = \emptyset$.

Then for $H^{ \delta}_{i,i+1}(z,\bar \alpha):=
\exists x_1 \dots \exists x_n \big(\underset{1\leq k\leq n}{\bigwedge}G_{k}(x_k,\bar \alpha) \land H_i (z, \bar x)\big)$ we have
$H_{i,i+1}^{\delta} (M,\bar \alpha)\cap V_{p,\varphi}(\beta)\neq\emptyset$,
$\neg H_{i,i+1}^{\delta} (M,\bar \alpha)\cap V_{p,\varphi}(\beta)\neq\emptyset$
and $V_{p,\varphi}(\delta_i)\subset H^{ \delta}_{i,i+1} (M,\bar \alpha)< V_{p,\varphi}(\delta_{i+1})$.
Next denote the convex 1-$\bar \alpha$-formula
$G_{i,i+1}^{\beta}(z,\bar \alpha):=\exists z_1\big(H_{i,i+1}^{ \delta}(z,\bar \alpha)\land \varphi (z_1,z)
\land \neg H_{i,i+1}^{\delta}(z_1,\bar \alpha)\big)$.

Thus, we obtain $n+1$ neighborhoods such that first (last) $n$ neighborhoods satisfy the initial condition.
\end{claim}

We continue this process of obtaining a new neighborhood $n$ times by Claim~\ref{CSS}:
removing first (last) neighborhood from $(n+1)$ neighborhoods, and considering the remaining $n$ as the initial.
Having obtained $2n$ neighborhoods, we divide them in half, and take the first $n$
as the next initial neighborhoods, and the second half -- as the second initial.
At the same time we define the convex 1-$\bar \alpha $-formula
$K_{\delta} (z,\bar \alpha):= G_1^{\delta} (M,\bar \alpha) < z < G_n^{\delta}(M, \bar \alpha)$.
After obtaining $2n$ neighborhoods, index the first $n$ neighborhoods by $\bar \gamma$ and the second $n$ neighborhoods
by $\bar \beta$.
Then define the following two convex 1-$\bar \alpha$-formulas:

$K_{\gamma} (z,\bar \alpha):= G_1^{\gamma} (M,\bar \alpha) < z < G_n^{\gamma}(M, \bar \alpha)$  and

$K_{\beta} (z,\bar \alpha):= G_1^{\beta} (M,\bar \alpha) < z < G_n^{\beta}(M, \bar \alpha)$.

Thus, $K_{\gamma} (M,\bar \alpha) < K_{\beta} (M,\bar \alpha)$.
Consequently, $K_{\gamma} (M,\bar \alpha) \cap K_{\beta} (M,\bar \alpha)\neq\emptyset$ and
\begin{center}
$(K_{\gamma} (M,\bar \alpha)\cup K_{\beta} (z,\bar \alpha)) \subset K_{\delta} (M,\bar \alpha)$.
\end{center}

So, if we can continue the process of Claim~\ref{CSS} infinitely many times,
 we obtain countable number of convex 1-$\bar \alpha$-formulas that is a base of continuum 1-types over $\bar\alpha$.
This contradicts to the assumption that the theory is small.
So impossibility to use the Claim~\ref{CSS} infinitely many times means that there exists
 $\bar \delta$ with type $tp(\bar \delta)$ satisfying $V_p$-$n$-condition.

\hfill $\Box$ Lemma~\ref{MoreConvexTypesQuasi}

Note that if the small theory $T$ is binary, then Lemmas~\ref{ConvexTypesQuasi} and~\ref{MoreConvexTypesQuasi}
do not require a proof.
And if $T$ is weakly o-minimal, for every $n\geq 2$
 the types $q_0$ and $q_0^n$ are complete and unique, and the proofs are simplified.

Let $tp(\bar \delta)=q(\bar x)$ satisfy $V_p$-$n$-condition then for any $i (1\leq i <n)$, for any 
$1$-$\bar \delta$-formula $H_i(z,\bar \delta)$ such that
$\mathfrak M \models \delta_i < H_i(M,\bar \delta)< \delta_{i+1}$ and 
\begin{center}
$\mathfrak M\models \forall x \forall y \Big(\exists z H_i(z,x,y, \bar \delta_n^{i,i+1}) \to
\forall z \big(H_i(z,x,y,\bar \delta_n^{i,i+1})\leftrightarrow H_i(z,x,y,\bar \delta_n^{i,i+1})^{+-} \big)\Big)$
\end{center}
 we have one of the two cases:

1) There exists  $k_0<\omega$ such that
$\forall z \big(\neg \varphi^{k_0}(z,x)\to \neg H_i (z,x,y,\bar \delta_n^{i,i+1})\big)\in q(x,y, \bar \delta_n^{i,i+1})$, 
or, equivalently, for
$T_1[H_i](x,y, \bar \delta_n^{i,i+1}):=\forall x_1\Big(\varphi(x,x_1) \to  
\exists z\big(H_i(z,x,y,\bar \delta_n^{i,i+1})\land \neg H_i(z,x_1,y,\bar \delta_n^{i,i+1})\big)\Big)
\land \forall y_1\Big(\varphi (y_1,y)\to \forall z\big(\neg H_i(z,x,y,\bar \delta_n^{i,i+1}) \leftrightarrow
\neg H_i(z,x,y_1, \bar \delta_n^{i,i+1})\big)\Big)$,
$T_1 [H_i](x,y, \bar \delta_n^{i,i+1})\in q (x,y,\bar \delta_n^{i,i+1})$.

2) There is  $l_0< \omega$ with
$\forall z \Big(\big(x< z \land \neg \varphi^{l_0}(y,z)\big)\to H_i(z,x,y, \bar \delta_n^{i,i+1})\Big)\in q(x,y,\bar \delta_n^{i,i+1})$, or, equivalently, for
$T_2 [H_i](x,y,  \bar \delta_n^{i,i+1}):=\forall y_1\Big(\varphi(y_1,y) \to\exists z\big(H_i(z,x,y, \bar \delta_n^{i,i+1})\land
\neg H_i(z,x,y_1, \bar \delta_n^{i,i+1})\big)\Big)\land
\forall x_1\Big(\varphi (x,x_1)\to \forall z\big(H_i(z,x,y, \bar \delta_n^{i,i+1})
\leftrightarrow H(z,x_1,y)\big)\Big)$, $T_2 [H_i](x,y, \bar \delta_n^{i,i+1})\in q(x,y, \bar \delta_n^{i,i+1})$.

Since 
\begin{center}
$\mathfrak M \models \forall x \forall y \big(T_1[H_i](x,y,  \bar \delta_n^{i,i+1}) \to \neg T_2 [H_i](x,y, \bar \delta_n^{i,i+1})\big) \land
\forall x \forall y\big(T_2[H_i](x,y, \bar \delta_n^{i,i+1}) \to\neg T_1 [H_i](x,y, \bar \delta_n^{i,i+1})\big)$,
\end{center}
we have that
$\big(T_1[H_i](x,y, \bar \delta_n^{i,i+1}) \lor T_2 [H_i](x,y, \bar \delta_n^{i,i+1})\big)\in q(x,y,\bar \delta_n^{i,i+1}).$

Therefore, the set of $n$-formulas
\begin{center}
$q_1^n(\bar x):= q_0^n(\bar x)\cup \Big\{T_1[H_i](\bar x) \lor T_2 [H_i](\bar x)\ |\
\mathfrak M \models \forall x \forall y \Big(H_i (x,x, y,\bar x_n^{i,i+1})\land H_i(M,\bar x) < y \land \forall z\big(H(z,x,y, \bar x_n^{i,i+1})\leftrightarrow H(z,x,y, \bar x_n^{i,i+1})^{+-}\big)\Big) \Big\}$
\end{center}
is consistent, and every complete extension of $q_1^n$ satisfies  $V_{p,\varphi}$-$n$-condition.
In fact, the proof of Lemma~\ref{MoreConvexTypesQuasi} is consistency of $q_1(x,y)$.

Let $S_n:=\{\langle i_1,i_2,\dots, i_n\rangle\ |\  i_1<i_2<\dots <i_n < \omega\}$, then  for any $\mu \in S_n$ denote
by $\bar x_{\mu}:=x_{\mu(1)} x_{\mu (2)}\dots x_{\mu (n)}$.
Consider $\Gamma (x_1, x_2, \dots, x_n, \dots)_{n<\omega}:= \underset{n<\omega,\ \mu \in S_n}{\bigcup} q_1^n (\bar x_{\mu}) $.

Consistence of $\Gamma$ follows from Lemma~\ref{MoreConvexTypesQuasi}.
Then there is a countable ordered set $B \subset N$ that satisfies $\Gamma$ 
and is ordered by the type of $\omega$, where $\mathfrak N$ is
an $\aleph_1$-saturated elementary extension of $\mathfrak M$.
Consider the Ehrefeucht-Mostovski type 
$EM(\omega/B)=\{\varphi(x_1,...,x_n)\ |$ for every $\mu\in S_n, n<\omega, \mathfrak N\models\varphi(\bar a_\mu)\}$.
Then, following the definition, $\Gamma(x_1, x_2, \dots, x_n, \dots)\subseteq EM(\omega/B)$. 
By the Standard  Lemma (\cite{TZ2012}, Lemma~5.1.3),
for any arbitrary infinite linear ordering there exists an indiscernible sequence  $\langle b_j\rangle_{j\in J}$.
Since $q_1^n (x_{\mu})$ is a subset of the Ehrefeucht-Mostovski type $EM (\omega/B)$,
every finite sequence of $J$ of length $n$ forms the tuple $\bar b_{\mu}$ that satisfies $V_p$-$n$-condition, 
what helps to implement the construction of continuum countable models.  

In the case when $T$ is weakly o-minimal, the ordered set $B$ is an indiscernible sequence
since for every $n\geq 2$ $q_0^n=q_1^n$ is a complete type.
When $T$ is binary, every indiscernible sequence of realizations of $\underset{n\geq 2}{\bigcup}q_0^n$ 
satisfies the $V_{p,\varphi}$-$n$-condition for every $n$.

Let $\tau:=\langle\tau_1, \tau_2, ..., \tau_i,...\rangle_{i<\omega}$,
$\tau_i\in\{0,1\}$, be an infinite sequence of zeros and ones.
Let $B_\tau:=\{\beta_1, \beta_2\}\cup
\{\beta_{2i-1,j}\ |\ i\in\mathbb N,j\in \mathbb Q\}\cup
\{\beta_{2i,1}, \beta_{2i,2}\ |\ i\in \mathbb N, \tau_i=0\}\cup$
$\{\beta_{2i,1}, \beta_{2i,2}, \beta_{2i,3}\ |\ i\in \mathbb N, \tau_i=1\}$
be an ordered indiscernible subset of $p(\mathfrak N)$ that exists by the previous statements.
The sets $V_{p(\mathfrak N)}(\beta_{i,j})$
are disjoint and ordered lexicographically by the indices $i,j$,
and $\beta_1<\beta_{i,j}<\beta_2$ for all $i$ and $j$.
Since the set $B_\tau$ is countable, fix an enumeration $B_\tau=\{b_1, b_2, ..., b_i,...\}$,
and construct a model $\mathfrak M_\tau:=\mathfrak A^{B_\tau}$ as in Theorem~\ref{Construction}.
By Lemma~\ref{MoreConvexTypesQuasi}, 
$p(\mathfrak M_\tau)\setminus\underset{\beta\in B_\tau}{\bigcup}V_{p(\mathfrak M)}(\beta)=\emptyset$.

Since the number of different infinite sequences $\tau$ of zeros and ones equals to
$2^{\aleph_0}$, and for every $\tau_1\not =\tau_2$,
$\mathfrak M_{\tau_1}\not\cong\mathfrak M_{\tau_2}$,
we have that $I(T\cup tp(\bar a\beta_1\beta_2),\aleph_0)=2^{\aleph_0}$.
By Proposition~\ref{Ext}, $I(T,\aleph_0)=2^{\aleph_0}$.

\hfill $\Box$ Theorem~\ref{abzQuasi}

\section{Extremely trivial types}\label{Strivial}

\begin{example}\label{Ex1}
Let   $\mathfrak{M}=\langle M;\ \Sigma\rangle$ be a countable structure of a finite signature $\Sigma$.
Let $C=\underset{a\in M}{\bigcup} C_a \subset \mathbb R$ be a countable set such that for every $a\not= b \in M$,
$C_a\cap C_b=\emptyset$ and $C_a$ and $C_b$ are mutually dense.
Consider $\mathfrak M^*=\langle C;\ \Sigma \cup \{<, E\}\rangle$, where
$E$ is a binary equivalence relation with $E(M^*,\ a)=C_a$ for every $a\in M$.
For an $n$-ary predicate $P\in \Sigma$, $a_1,a_2,...,a_n \in M$, and
$b_1\in C_{a_1}, b_2\in C_{a_1},  ..., b_n\in C_{a_n}$
define $\mathfrak M^*\models P(b_1,b_2,...,b_n)$ if
 $\mathfrak{M}\models P(a_1,a_2,...,a_n)$.
Functions and constants are defined analogically.
It is easy to see that
$\langle C/E;\ \Sigma \rangle\cong \langle M;\ \Sigma \rangle$, and then
$I(Th(\mathfrak M),\aleph_0)=I(Th(\mathfrak M^*),\aleph_0)$.
By this, for an arbitrary structure, we can define a structure with an artificial linear order that
has the same number of countable non-isomorphic models.
\end{example}

\begin{definition}\label{ExtremelyTrivial}\cite{BBZ2018}
Let $T$ be a small complete theory,
 $p(\bar x)$ be a non-principal type over a finite subset $A$ of some model of $T$.
The type $p$ is {\bf extremely trivial}, if for every natural number $n\geq 1$ and every $n$ realizations
$\bar\beta_1, \bar\beta_2$, ...$\bar\beta_n$ of $p$,
$p(\mathfrak M(\bar\beta_1, \bar\beta_2,...,\bar\beta_n,\bar a))
=\{\bar\beta_1, \bar\beta_2...\bar\beta_n\}$, where
$\bar a$ is some enumeration of the set $A$,
and $\mathfrak M(\bar\beta_1, \bar\beta_2,...,\bar\beta_n,\bar a)$
is a prime model over $\{\bar\beta_1, \bar\beta_2,...,\bar\beta_n,\bar a\}$.
\end{definition}

\begin{example}\label{Ex2}
Let $T$ be a linearly ordered theory, $\mathfrak M\models T$ be sufficiently saturated, and $A\subset M$ be finite.
Let $\theta(x)$ be a convex $A$-formula, and $\{\varphi_i(x)\ | i<\omega\}$ be a countable set of 1-$A$-formulas such that 

1) for every $i$, $\varphi_i(M)$ is mutually dense with $\neg\varphi_i(M)$ in $\theta(M)$,

2) for every $i,j<\omega$, $i\neq j$, $\varphi_i(M)$ and $\varphi_j(M)$ are mutually dense in $\theta(M)$, and

3) there is no infinite family of uniformly definable equivalence relations $\varepsilon(x,y,\bar b_i)$ ($i<\omega$, $\bar b_i\in M$) 
on $\theta(M)$, such that for every $i<\omega$ every
the set of all $\varepsilon(x,y,\bar b_i)$-classes is infinite and densely ordered, 
and for every $i\neq j$ all $\varepsilon(x,y,\bar b_i)$- and $\varepsilon(x,y,\bar b_j)$-classes are disjoint.

Suppose that $p(x):=\theta(x)\cup\{\neg \varphi_i(x)\ |\ i<\omega\}$ has no convex to the right and 
to the left formulas.
Then every completion of $p$ to a type from $S_1(A)$ is extremely trivial.
\end{example}

\begin{theorem}\label{Trivial}\cite{BBZ2018}
Let $T$ be a countable complete linearly ordered theory.
If there exists a finite subset  $A$ of a model $\mathfrak M\models T$
and a non-principal extremely trivial type $p(x)\in S_1(A)$,
then $T$ has $2^{\aleph_0}$ countable non-isomorphic models.
\end{theorem}

Theorem~\ref{Trivial} was presented in~\cite{BBZ2018}, and is one of a series of results obtained on the way to the current state of
construction in Theorem~\ref{Construction}.
Any realization of an extremely trivial type has a non-principal type over any other its realizations.
And construction in Theorem~\ref{Construction}, while adding only realizations of principal formulas,
guarantees omitting those non-principal types.
By this, we are allowed to arrange the realization set of the extremely trivial type in $2^{\aleph_0}$ non-isomorphic ways.

\begin{example}\label{Ex3}
Let $\mathfrak M=\langle M; =, <, \varepsilon, c_i\rangle_{i<\omega}$, where
$<$ is a binary relation of linear order, $c_i$ are constants with $c_i<c_{i+1}$,
$\varepsilon$ is a binary equivalence relation such that
every $\varepsilon$-class is infinite and dense,
and every two $\varepsilon$-classes are disjoint and mutually dense.
In $Th(\mathfrak M)$ there is a countable number of non-extremely-trivial non-isolated types:
for each $n<\omega$, $p_n:=\{x>c_i\}_{i<\omega}\cup\{\varepsilon(x,c_n)\}$,
and  $p_0:=\{x>c_i\}_{i<\omega}\cup\{\neg\varepsilon(x,c_i)\}_{i<\omega}$.

We can extend the signature of $\mathfrak M$ to include relations of any complexity
that respect equivalence classes of $\varepsilon$.
Then, for a sufficiently saturated elementary extension $\mathfrak N$ of $\mathfrak M$,
$p^c_0(\mathfrak N)/\varepsilon$ defines a structure, in which the linear order is not of any importance.
\end{example}

Examples~\ref{Ex1} and~\ref{Ex3} illustrate that studying linear orders is equivalent to studying
model theory in general.
Although these cases are promising for research,
at this stage we restrict ourselves to theories with no
uniformly definable infinite families of mutually dense sets.

Under this restriction, if some 1-type over a finite set has an infinite number of types with the same convex closure,
then, by Theorem~\ref{Trivial}, an extremely trivial type appears, what leads
to the maximal number of countable models.
Later we consider the case, when the number of types with the same convex closure can not be infinite.

\section{Equivalence-generating formulas}\label{SEG}

As before, we work with a countable complete linearly ordered theory $T$.

\begin{definition}\cite{BK2006, Baizh2001}
Let $\mathfrak M$ be a linearly ordered structure,
$A\subseteq M$, $M$ be $|A|^+$-saturated, and $p\in S_1(A)$ be non-algebraic.
A $p$-preserving convex to the right (left) $A$-definable formula $\varphi(x,y)$ is
{\bf equivalence-generating} if for every $\alpha\in p(M)$ and every $\beta\in \varphi(M,\alpha)\cap p(M)$ the following holds:
\begin{center}
$\mathfrak M\models \forall x (x\geq\beta\rightarrow(\varphi(x,\alpha)\leftrightarrow\varphi(x,\beta)))$

($\mathfrak M\models \forall x (x\leq\beta\rightarrow(\varphi(x,\alpha)\leftrightarrow\varphi(x,\beta)))$).
\end{center}
\end{definition}

By $CRF(p)$ ($CLF(p)$) we denote the family of all $p$-preserving convex to the right
(left) $A$-formulas.
In weakly o-minimal theories of finite convexity rank $CRF(p)$ and $CLF(p)$ are always finite.
We say that $\varphi_1(x,y)\in CRF(p)$ is {\bf greater} than $\varphi_2(x,y)\in CRF(p)$ on $p$, if
$\varphi_2(M,\alpha)\cap p(\mathfrak M)\subset \varphi_1(M,\alpha)\cap p(\mathfrak M)$ for some
(equivalently, for every) $\alpha\in p(\mathfrak M)$.
Formulas with $\varphi_2(M,\alpha)\cap p(\mathfrak M)= \varphi_1(M,\alpha)\cap p(\mathfrak M)$
 we consider equivalent on $p$.

\smallskip

\noindent{\bf Restriction}. {\it
We restrict to a linearly ordered theory $T$ such that
for every subsets
$A$ and $B$ of a model of $T$, where
$A$ is finite and $B$ is $A$-definable

1)  for every $A$-formula
$E(x,y,\bar z)$
there is no infinite sequence
$\bar b_1, \bar b_2, ..., \bar b_i,... \in B$ such that for every $i<\omega$,
$E(x,y,\bar b_i)$ is an equivalence relation on $p$ with
convex in $p$ mutually dense classes partitioned into infinitely many infinite
$E(x,y,\bar b_{i+1})$-classes;

2) the set $\{q\in S_1(A) \ |\ p^c=q^c\}$ is finite for every 1-type $p\in S_1(A)$;

3) for every $p\in S_1(A)$ all $p$-preserving convex to the right
and all $p$-preserving convex to the left
formulas are equivalence-generating.}

\begin{theorem}\label{NoGreatest}\cite{BUZ2019b}
Let $T$ be a countable complete linearly ordered theory satisfying the Restriction,
$A$ be a finite subset of a countable saturated model $\mathfrak M$ of $T$,
and $p(x)\in S_1(A)$ be a non-algebraic 1-type.
If $CRF(p)$ is infinite and has no greatest formula, then
$T$ has $2^{\aleph_0}$ countable non-isomorphic models.
\end{theorem}

\noindent{}{\it Proof of Theorem~\ref{NoGreatest}.}
Since every non-small theory has $2^{\aleph_0}$ countable non-isomorphic models,
it remains to prove the case, when the theory $T$ is small.
For simplicity we extend our language to $\mathcal L(A)$ and work in the theory
$T\cup tp(\bar a)$, where $\bar a$ is an enumeration of the set $A$.

For an arbitrary $p$-preserving  convex  to the right $2$-formula
we define an equivalence relation such that any its class is convex in $p(\mathfrak M)$.
For $\varphi(x,y)\in CRF(p)$ denote by $\varepsilon_{\varphi}(x, y)$ the formula
$\big(\varphi (y,x)\lor\varphi(x,y)\big)$.
It can be proved that this formula is equivalent to
$\forall z(\varphi (z,x)^+\leftrightarrow\varphi(z, y)^+)$.
The formula $\varepsilon_{\varphi}(x,y)$
defines an equivalence relation with convex classes on $p(\mathfrak M)$,
but not necessarily on $p^c(\mathfrak M)$.

Since  $|\{q\in S_1(A)\ |\ p^c=q^c\}|<\omega$, then for some $\Theta \in p$,
$p(\mathfrak M)= p^c(\mathfrak M)\cap \Theta(M)$.
Then, $\varepsilon_{\varphi}$ is an $A$-definable equivalence relation with convex
$\varepsilon_{\varphi}$-classes on $\Theta(M)$ for suitable $\Theta \in p$.
Thus, for arbitrary $\alpha \in p(\mathfrak M)$, $\varepsilon_{\varphi}(M,\alpha)\cap \Theta(M) $
is exactly $A\cup\{\alpha\}$-definable convex $\varepsilon_{\varphi}$-class on
$p(\mathfrak M)$ containing $\alpha$.

Denote $\varepsilon_{\varphi, \Theta} (M,\alpha)=\varepsilon_{\varphi}(M,\alpha)\cap \Theta(M)$.
Let $\varphi_1,\varphi_2\in CRF(p)$.
It is easy to see that for some (equivalently, for every) $\alpha\in p(\mathfrak M)$
 $\varphi_1(M,\alpha)\subseteq\varphi_2(M,\alpha)$ if and only if
$\varepsilon_{\varphi_1,\Theta}(\alpha')\subseteq\varepsilon_{\varphi_2,\Theta}(\alpha')$
for every $\alpha'\in p(\mathfrak M)$.

Note that it is possible to define a corresponding equivalence relation on $p^c(\mathfrak M)$ by letting
$\varepsilon^1_{\varphi}(x, y):= \exists z\bigg(\Theta(z)\land\forall z_1\forall z_2\Big(
\big(\varepsilon_{\varphi, \Theta}(z_1,z)^-\land\varepsilon_{\varphi, \Theta}(z_2,z)^+\big)\rightarrow
\big(z_1<x<z_2\land z_1<y<z_2\big)\Big)\bigg)$.
This as well is an $\emptyset$-definable equivalence relation with convex classes on some formula from $p^c$.

For $\alpha\in p(\mathfrak M)$ let
$V_{p(\mathfrak M)}(\alpha):=\{\gamma\in p(\mathfrak M)\ |$ there exists $\varphi(x,y)\in CRF(p)$ such that
$\mathfrak M\models \varphi(\alpha,\gamma)\lor\varphi(\gamma,\alpha)\}$.
It follows from the definition that
\begin{center}
$V_{p(\mathfrak M)}(\alpha)=\underset{\varphi \in CRF(p)}{\bigcup}\varepsilon_{\varphi, \Theta}(M,\alpha)$.
\end{center}

Since there is no greatest equivalence-generating formula, the sets
$V_{p(\mathfrak M)}(\alpha)$, $V_{p(\mathfrak M)}(\alpha)^+$ and
$V_{p(\mathfrak M)}(\alpha)^-$ are not $A\alpha$-definable.

For arbitrary $\alpha \in p^c (\mathfrak M)$  denote
$V_{p^c(\mathfrak M)}(\alpha):=\cup_{\varphi \in CRF(p)}\varepsilon^1_{\varphi}(M,\alpha)$.
It follows from the definition that for any $\alpha \in p(\mathfrak M)$

\begin{center}
$V_{p(\mathfrak M)}(\alpha)\subseteq V_{p^c(\mathfrak M)}(\alpha),\;
V_{p(\mathfrak M)}(\alpha)^+= V_{p^c(\mathfrak M)}(\alpha)^+, \;
V_{p(\mathfrak M)}(\alpha)^-=V_{p^c(\mathfrak M)}(\alpha)^-$.
\end{center}

Denote
\begin{center}
$(V_{p}(\alpha),V_{p}(\beta))_{p(\mathfrak M)}:=
\{\gamma\in p(\mathfrak M) \ |\ V _{p(\mathfrak M)}(\alpha)<\gamma<V_{p(\mathfrak M)}(\beta)\}$.
\end{center}
Then
\begin{center}
$(V_{p}(\alpha),V_{p}(\beta))_{p(\mathfrak M)} =\{\gamma \in \Theta (M)\ |$
for every $\varphi\in CRF(p)$,
$\mathfrak M \models \alpha <\gamma< \beta \land\neg \varepsilon_{\varphi}(\alpha,\gamma)\land
\neg \varepsilon_{\varphi}(\beta,\gamma)\}$.
\end{center}

\begin{lemma}\label{ConvexTypes}
Let $\alpha, \beta\in p(\mathfrak M)$ be such that $V_{p(\mathfrak M)}(\alpha)<V_{p(\mathfrak M)}(\beta)$.
Then for all realizations $\gamma_1,\gamma_2\in (V_{p}(\alpha),V_{p}(\beta))_{p(\mathfrak M)}$
$tp^c(\gamma_1/\alpha\beta)=tp^c(\gamma_2/\alpha\beta)$.

\end{lemma}

\noindent{}\noindent{\it Proof of Lemma~\ref{ConvexTypes}}.
Assume that the conclusion of the lemma is not true.
This means there is a $\{\alpha, \beta\}$-definable convex $1$-formula $H(x,\alpha,\beta)$
that distinguishes $\gamma_1$ and $\gamma_2$.
We can arrange that $\alpha,\gamma_1\in H(M,\alpha,\beta)<\gamma_2$ and
$\alpha$ is the left endpoint of $H(M,\alpha,\beta)$.

Let $p_1(x, \beta):=tp(\alpha/\beta)$.
It follows from the definition that
$p(x)\cup (x< V_{p(\mathfrak M)}(\beta))=
p(x)\cup \{x<\varepsilon_{\varphi,\Theta}(M,\beta)\ |\ \varphi \in CRF(p)\}\subseteq p_1$
and $p_1^c(\mathfrak M)=p^c(\mathfrak M)\cap V_{p(\mathfrak M)}(\beta)^-$.
By the second condition of the Restriction the set $\{q_1\in S_1(\{\beta\})\ |\ p_1^c=q_1^c\}$ is finite.
This means there is a $\{\beta\}$-definable zebra-formula $\Theta_1\in p_1$ such that $p_1^c\cup\{\Theta_1\}= p_1$.
We have that $H(x, y,\beta)$ is a convex to the right $p_1$-preserving $\{\beta\}$-formula.
Notice that $CRF(p)\cup \{H(x,y,\beta)\}\subseteq CRF(p_1)$
since for every $\varphi \in CRF (p)$, $H(M, \alpha, \beta)^+\subset \varphi (N,\alpha)^+$.

Then by the third condition of the Restriction the $\{\beta\}$-definable formula $H(x, y, \beta)$
defines an equivalence relation $\varepsilon_{H(x, y, \beta)}(x,y, \beta)$ with
convex $\varepsilon_H$-classes on the set of all realizations of $p_1$.
Define a new equivalence relation on some convex part of $p(\mathfrak M)$:
\begin{center}
$E_{H, \Theta}(x,y, \beta):= \exists z \big(\Theta_1(z)\land \Theta (x)\land \Theta (y)
\land z\leq x \land z\leq y \land \neg H(x, z, \beta)^+ \land \neg H(y, z, \beta)^+\big)$.
\end{center}

Last sentence gives us a $\{\beta\}$-definable equivalence relation on $V_{p(\mathfrak N)}^-(\beta)$.
Since $\gamma_1\not \in H(N, \alpha, \beta)^+$,
$\mathfrak N \models E_{H, \Theta}(\gamma_1, \alpha, \beta)$.
Since for $\varphi \in CRF(p)$,  $\mathfrak N\models \neg \varepsilon_{\varphi, \Theta}(\gamma_1, \alpha)$,
the definable set $E_{H, \Theta}(M,\alpha, \beta)$ is a convex subset of
$p(\mathfrak M)$, that contains a densely ordered infinite set of non-definable sets of $V_{p(\mathfrak M)}$-kind.

Consider the set of $1$-formulas $p_1(x, \alpha)$. This is a complete 1-type because $\alpha, \beta \in p(\mathfrak M)$.
Take an arbitrary realization of $p_1(x, \alpha)$  outside of  $ E_{H,\Theta}(M, \alpha, \beta)$ and denote it by $\alpha_1$.
Then either we have $E_{H, \Theta}(M, \alpha_1, \alpha)\subset E_{H, \Theta}(M, \alpha, \beta)$,
or $E_{H, \Theta}(M, \alpha_1, \alpha)\subset E_{H, \Theta}(M, \alpha_1, \beta)$.
Suppose that the first is true.
The  set $E_{H, \Theta}(M, \alpha, \beta)$, that is convex in $p(\mathfrak M)$,
contains an infinite number of convex sets definable by $E_{H,\Theta}(x, y,\alpha)$.

Let $f$ be an $A$-isomorphism of the countable saturated model $\mathfrak N$
 generated by an elementary monomorphism such that $f(\beta)=\alpha$ and $f(\alpha)=\alpha_1$.
The existence of such an elementary monomorphism follows from $tp(\alpha, \beta)=tp(\alpha_1, \alpha)$.
The isomorphism $f $
generates an infinite sequence $\langle\alpha_n\rangle_{0<n<\omega}$ such that for every $n<\omega$, we have
\begin{center}
$E_{H, \Theta}(M, \alpha_{n+2}, \alpha_{n+1})\subset E_{H, \Theta}(M, \alpha_{n+1}, \alpha_{n})$.
\end{center}
This contradicts to the first condition of the Restriction.

$\hfill\Box$ Lemma~\ref{ConvexTypes}

By Lemma~\ref{ConvexTypes} the type $tp^c(\gamma_1/\alpha \beta)$ is irrational,
and therefore $tp(\gamma_1/\alpha \beta)$ is non-principal.

A corollary of the proof of Lemma~\ref{ConvexTypes}
is the following lemma.

\begin{lemma}\label{MoreConvexTypes}
Let $n<\omega$, $\alpha_i\in p(\mathfrak M)$, $i\leq n$, such that
$V_{p(\mathfrak M)}(\alpha_i)<V_{p(\mathfrak M)}(\alpha_{i+1})$ for all $1\leq i\leq n-1$.
If $\gamma_1,\gamma_2\in p(\mathfrak M)$ are such that for every $i$, $1\leq i\leq n$,
$\gamma_1,\gamma_2\not\in V_{p(\mathfrak M)}(\alpha_i)$, and
$\gamma_1<\alpha_i$ if and only if $\gamma_2<\alpha_i$,
then $tp^c(\gamma_1/\alpha_1,...,\alpha_n)=
tp^c(\gamma_2/\alpha_1,...,\alpha_n)$.

\end{lemma}

Lemma~\ref{MoreConvexTypes} implies that the type $tp(\gamma_1/A\cup\{\alpha_1,...,\alpha_n\})$ is non-principal.

Let $\tau:=\langle\tau_1, \tau_2, ..., \tau_i,...\rangle_{i<\omega}$,
$\tau_i\in\{0,1\}$, be an infinite sequence of zeroes and ones.

Let $\mathfrak N$ be an $\aleph_1$-saturated elementary extension of $M$.
Let $B_\tau:=\{\beta_1,\beta_2\}\cup\{\beta_{2n-1,q}\ |\ n\in\mathbb N,q\in \mathbb Q\}\cup
\{\beta_{2n,1}, \beta_{2n,2}\ |\ n\in \mathbb N, \tau_n=0\}\cup$
$\{\beta_{2n,1}, \beta_{2n,2}, \beta_{2n,3}\ |\ n\in \mathbb N, \tau_n=1\}$
be a subset of $p(\mathfrak N)$ such that the sets $V_{p(\mathfrak N)}(\beta_{n,q})$
are disjoint and ordered lexicographically by the indices $n,q$,
and $V_{p(\mathfrak N)}(\beta_1)<V_{p(\mathfrak N)}(\beta_{n,q})<V_{p(\mathfrak N)}(\beta_2)$ for every $\beta_{n,q}$.
So, $B_\tau$ consists of densely ordered neighborhoods $V_{p(\mathfrak N)}$ with inclusions of discrete chains
of neighborhoods of lengths depending on $\tau$, all of them bounded by $V_{p(\mathfrak N)}(\beta_1)$ and $V_{p(\mathfrak N)}(\beta_2)$.
Apply Theorem~\ref{Construction} to construct a model $\mathfrak M_\tau:=\mathfrak A^{B_\tau}$.

Denote by $p'\in S_1(\{\beta_1,\beta_2\})$ some compete type extending $p(x)\cup\{\beta_1<x<\beta_2\}$.
We claim that $p'(\mathfrak M_\tau)\setminus\underset{b\in B_\tau}
{\bigcup}V_{p'(\mathfrak M_\tau)}(b)=\emptyset$.
Towards a contradiction suppose that there exists a realization
$\gamma\in p'(\mathfrak M_\tau)\setminus\underset{\beta\in B_\tau}
{\bigcup}V_{p'(\mathfrak M_\tau)}(\beta)$.
Then $\gamma\in(V_p(\beta_1),V_p(\beta_2))_{p(\mathfrak M_\tau)}$.
By Theorem~\ref{Construction}, 3), the type $tp(\gamma/\bar \beta)$ is
principal for some $\bar \beta\in B_\tau$.
By constructions in Lemma~\ref{chain} and Theorem~\ref{Construction} there is
$\bar\beta'\in B_\tau$, that includes $\beta_1$, $\beta_2$ and $\bar \beta$, and such that
$tp(\gamma/\bar \beta')$ is also principal.
But by Lemma~\ref{MoreConvexTypes}
the type $tp(\gamma/\bar \beta')$ is non-principal.
This is a contradiction and we have that
$p'(\mathfrak M_\tau)\setminus\underset{b\in B_\tau}{\bigcup}V_{p'(\mathfrak M_\tau)}(\beta)=\emptyset$.

Then, the order type of neighborhoods in $p'(\mathfrak M_\tau)$ is determined by the order type of $B_\tau$,
and models based on different infinite sequences $\tau$ are non-isomorphic.
Since the number of different sequences $\tau$ equals to $2^{\aleph_0}$,
$I(T\cup tp(\bar a, \beta_1,\beta_2),\aleph_0)=2^{\aleph_0}$,
and, by Proposition~\ref{Ext}, $I(T,\aleph_0)=2^{\aleph_0}$.

$\hfill\Box$ Theorem~\ref{NoGreatest}

\begin{theorem}\label{NoLeast}
Let $T$ be a countable complete linearly ordered theory satisfying the Restriction,
$A$ be a finite subset of a countable saturated model $\mathfrak M$ of $T$,
and $p(x)\in S_1(A)$ be a non-algebraic 1-type.
If $CRF(p)$ is infinite and has no least formula, then
$T$ has $2^{\aleph_0}$ countable non-isomorphic models.
\end{theorem}

\noindent{\it Proof of Theorem~\ref{NoLeast}.}
As in Theorem~\ref{NoGreatest}, we consider $T$ to be small, and extend its language to $\mathcal L(A)$.
Let $\bar a$ be a tuple enumerating $A$.
Let $\Theta \in p$ be such that $p(\mathfrak M)= p^c(M)\cap \Theta(M)$.
We use notations introduced in the proof of Theorem~\ref{NoGreatest}.

For $\alpha\in p(\mathfrak M)$ define
 $Ker_{p(\mathfrak M)}(\alpha):=\underset{\varphi\in CRF(p)}{\bigcap}\varepsilon_{\varphi, \Theta}(N,\alpha)$
to be its kernel in $p(\mathfrak M)$.
It is obvious that the relation $\alpha\in Ker_{p(\mathfrak M)}(\beta)$ is an equivalence relation
with convex classes on $p(\mathfrak M)$.

\begin{lemma}\label{OneKernel}
Let $\alpha_1,\alpha_2\in p(\mathfrak M)$,
$\alpha_1\neq\alpha_2$
and $\alpha_2\in Ker_{p(\mathfrak M)}(\alpha_1)$.
Then $tp(\alpha_2/\alpha_1)$ is non-principal.
\end{lemma}

\noindent{\it Proof of Lemma~\ref{OneKernel}.}

Notice that since all convex to the right and convex to the left $p$-preserving formulas are
equivalence-generating, there is no least formula in $CLF(p)$ as well.
We suppose that $\alpha_1<\alpha_2$. 
The case $\alpha_1>\alpha_2$ can be considered analogically.
Towards a contradiction suppose that $H(x,\alpha_1)$ is an isolating formula of $tp(\alpha_2/\alpha_1)$.
For every $\varphi\in CRF(p)$ $\varepsilon_{\varphi,\Theta}(x,\alpha_1)\in tp(\alpha_2/\alpha_1)$,
therefore $H(M,\alpha_1)\subseteq \varepsilon_{\varphi,\Theta}(M,\alpha_1)$.
Then $H^c(M,\alpha_1)$ as well as the definable set of
$(H^c(x,\alpha_1)\land x\geq\alpha_1)$
is a subset of $Ker_{p(\mathfrak M)}(\alpha_1)$.
But since the last formula is convex to the right $p$, 
and $CRF(p)$ and haves no least formula, this is a contradiction.

$\hfill\Box$ Lemma~\ref{OneKernel}

\begin{lemma}\label{ManyKernels}
Let $n,m<\omega$, $\alpha_1,...,\alpha_n, \alpha_{n+1}$ be an increasing sequence from $p(\mathfrak M)$,
 $\beta_1,...,\beta_m\in p(\mathfrak M)$,
and for $i\leq n$, $j\leq m$
$\alpha_i\in Ker_{p(\mathfrak M)}(\alpha_{n+1})$ and $\beta_j\not\in Ker_{p(\mathfrak M)}(\alpha_{n+1})$.
Then the types $tp(\alpha_{1}/\alpha_2 \alpha_2...\alpha_n\bar\beta_n)$ and 
$tp(\alpha_{n+1}/\bar\alpha_n\bar\beta_n)$ are non-principal, where
$\bar\alpha_n:=(\alpha_1,...,\alpha_n)$ and $\bar\beta_n:=(\beta_1,...,\beta_m)$.
\end{lemma}

\noindent{\it Proof of Lemma~\ref{ManyKernels}.}
We give a proof for $tp(\alpha_{n+1}/\bar\alpha_n\bar\beta_n)$, the proof
for $tp(\alpha_{1}/\alpha_2 \alpha_2...\alpha_n\bar\beta_n)$ is similar.

The proof is by induction.
The first step is done by Lemma~\ref{OneKernel}.

It is easy to show that the type $tp(\alpha_{n+1}/\bar\alpha_n)$ is non-principal, as well as that
$tp(\alpha_{n+1}/\bar\alpha_n\bar\beta_n)$ is non-principal in the case when 
all elements of $\bar\beta_n$ are all in different kernels:
the type $tp(\bar\beta_n/\bar\alpha_n)$ is principal and 
the type $tp(\alpha_{n+1}/\bar\alpha_n)$ is non-principal, 
but a set definable with parameters from a principal type can not 
be a subset of the set of all realizations of a non-principal type.

Without loss of generality suppose that all elements of the tuple $\bar\beta_n$
are in the same kernel.
Obviously, $tp(\alpha_{n+1}/\bar\alpha_n\beta_1)$ is non-principal.
Suppose that $tp(\alpha_{n+1}/\bar\alpha_n\bar\beta_{m})$ is principal.
And, towards a contradiction, suppose that $H(x,\bar\alpha_n,\bar\beta_{m+1})$
is an isolating formula of the type $tp(\alpha_{n+1}/\bar\alpha_n\bar\beta_{m+1})$.
Then $(H^c(x,\bar\alpha_n,\bar\beta_{m+1})\land x\geq\alpha_n)$ is convex to the right on 
$tp(\alpha_{n+1}/\bar\alpha_n\bar\beta_{m})$.
Let $\varepsilon(x,y,\bar\alpha_n,\bar\beta_{m+1})$ be the equivalence relation
 generated by the previous formula.
Define a new equivalence relation
\begin{center}
$R(x,y,\bar\alpha_n,\bar\beta_{m}):=\forall x_1\forall y_1\big(
\varepsilon(x_1,y_1,\bar\alpha_n,\bar\beta_{m},x)\leftrightarrow
\varepsilon(x_1,y_1,\bar\alpha_n,\bar\beta_{m},y)\big)$.
\end{center}

Take $\beta^1$ with $tp(\beta^1/\bar\alpha_n,\bar\beta_{m})=
tp(\beta_{m+1}/\bar\alpha_n,\bar\beta_{m})$
such that $\beta^1$ is in a different $\varepsilon$-class than every other $b_i$,
either every $\varepsilon(x,y,\bar\alpha_n,\bar\beta_{m},\beta^1)$-class is 
a subset, or a superset of $\varepsilon(x,y,\bar\alpha_n,\bar\beta_{m+1})$.
Next, take $\beta^2$ be in a different $\varepsilon$-class from all $\beta_i$'s and $\beta^1$.
We can repeat this procedure infinitely many times since all of the
$R(x,y,\bar\alpha_n,\bar\beta_{m})$-classes are convex and densely ordered.
Each next $b^i$ will either expand, or narrow down its $\varepsilon(x,y,\bar\alpha_n,\bar\beta_{m},\beta^i)$-classes.
This way, we obtain a contradiction with Restriction, part 1).

$\hfill\Box$ Lemma~\ref{ManyKernels}

Let $\mathfrak N$ be an $\aleph_1$-saturated elementary extension of $\mathfrak M$.
There are five kinds of kernels in $p(\mathfrak N)$: a singleton and infinite kernels with or without right or left border.
Fix $\tau:=\langle\tau_1, \tau_2, ..., \tau_i,...\rangle_{i<\omega}$,
$\tau_i\in\{0,1\}$, be an infinite sequence of zeroes and ones.
We construct a set $B_\tau$ with different kinds of kernels arranged depending on $\tau$.
Let $B_\tau:=\{\beta_{2n-1,q,1}\ |\ n\in\mathbb N, q\in\mathbb Q\}\cup
\{\beta_{2n,1,m}\ |\ n,m\in\mathbb N, \tau_n=0\}\cup
\{\beta_{2n,1,m}\ |\ n\in\mathbb N, m\in\mathbb N^- \tau_n=1\}$
be a subset of $p(\mathfrak N)$ such that $\beta_{n_1,q_1,m_1}\in Ker_{p(\mathfrak N)}(\beta_{n_2,q_2,m_2})$
if and only if $n_1=n_2$ and $q_1=q_2$,
all the $\beta_{n,q,m}$ are ordered lexicographically by their indices,
and for every element from $p(\mathfrak N)$ there is a representative of its kernel in $B_\tau$.
Apply Theorem~\ref{Construction} to construct a model $\mathfrak M_\tau:=\mathfrak A^{B_\tau}$.

Suppose  there exists $\gamma\in p(\mathfrak M_\tau)\setminus B_\tau$ such that
$\gamma\in Ker_{p(\mathfrak M_\tau)}(\alpha)$  and it is larger (smaller) than every element of
$B_\tau\cap Ker_{p(\mathfrak M_\tau)}(\alpha)$ for some $\alpha\in B_\tau$.
By Theorem~\ref{Construction}, 3), there is $\bar \beta\in B_\tau$ such that
$tp(\gamma/\bar \beta)$ is principal.
By the proofs of Lemma~\ref{chain} and Theorem~\ref{Construction}
there is $\bar \beta'\in B_\tau$ extending $\bar \beta$ and
such that $\bar \beta'$ includes at least one element form $Ker_{p(\mathfrak M_\tau)}(\alpha)\cap B_\tau$
and such that $tp(\gamma/\bar \beta')$ is also principal.
But by Lemma~\ref{ManyKernels} $tp(\gamma/\bar \beta')$ is non-principal.
This is a contradiction, and, therefore, there is no such $\gamma$.
Then kernels in $p(\mathfrak M_\tau)$ are of the same kinds as the corresponding kernels in $B_\tau$.

Now, suppose that a new kernel appears during construction.
That is, there exists $\gamma\in p(\mathfrak (M_\tau)\setminus \underset{\beta\in B_\tau}{\bigcup}Ker_{p(\mathfrak (M_\tau)}(\beta)$.
Then, by Lemma~\ref{ManyKernels}, $|Ker_{p(\mathfrak (M_\tau)}(\gamma)|=1$.
Therefore, appearance of $\gamma$ does not affect the order type:
the order types of kernels in $p(\mathfrak (M_\tau)$ are determined by their order types in $B_\tau$.
Then, $\mathfrak M_{\tau_1}\not\cong\mathfrak M_{\tau_2}$ whenever $\tau_1\neq\tau_2$.

Since the number of different sequences $\tau$ equals to
$2^{\aleph_0}$, $I(T\cup tp(\bar a),\aleph_0)=2^{\aleph_0}$.
By Proposition~\ref{Ext} we obtain $I(T,\aleph_0)=2^{\aleph_0}$.

$\hfill\Box$ Theorem~\ref{NoLeast}

The following theorem generalizes theorems~\ref{NoGreatest} and \ref{NoLeast}.

\begin{theorem}\label{Infinite}
Let $T$ be a countable complete linearly ordered theory satisfying the Restriction,
$A$ be a finite subset of a countable saturated model $\mathfrak M$ of $T$,
and $p(x)\in S_1(A)$ be a non-algebraic 1-type.
If $CRF(p)$ is infinite, then
$T$ has $2^{\aleph_0}$ countable non-isomorphic models.
\end{theorem}

\noindent{\it Proof of Theorem~\ref{Infinite}.}
As in previous theorems, we consider $T$ to be small, and extend its language to $\mathcal L(A)$.
Let $\bar a$ be a tuple enumerating $A$,
and $\Theta \in p$ be such that $p(\mathfrak M)= p^c(M)\cap \Theta(M)$.
For $\varphi(x,y)\in CRF(p)$ and $\alpha, \beta \in p(\mathfrak M)$ denote $\varepsilon_{\varphi,\Theta}(\alpha, \beta)$
as in Theorem~\ref{NoGreatest}.

Theorems~\ref{NoGreatest} and~\ref{NoLeast} deal with the cases, when $CRF(p)$ has
no greatest, or no least formula.
So, we suppose $CRF(p)$ has both a greatest formula, $\varphi_{max}$, and a least formula, $\varphi_{min}$.
We consider two cases.

{\bf Case 1.}
There is a formula $\psi(x,y)\in CRF(p)$ such that for every $\psi_1\in CRF(p)$ with 
$\psi_1(M,\alpha)\cap p(\mathfrak M)\subset \psi(M,\alpha)$
for some (equivalently, for every) $\alpha\in p(\mathfrak M)$
there is $\psi_2\in CRF(p)$ such that  $\psi_1(M,\alpha)\cap p(\mathfrak M)\subset\psi_2(M,\alpha)
\subset\psi(M,\alpha)\cap p(\mathfrak M)$.

For $\alpha,\beta\in p(\mathfrak M)$ denote
\begin{center}
$V^\psi_{p(\mathfrak M)}(\alpha):=\underset{\varphi(M,\alpha)\cap p(\mathfrak N)\subset\psi(M,\alpha)\cap p(\mathfrak M)}{\underset{\varphi\in CRF(p),}{\bigcup}}
\varepsilon_{\varphi,\Theta}(M,\alpha)$,

$(V^\psi_{p}(\alpha),V^\psi_{p}(\beta))_{p(\mathfrak M)}:=
\{\gamma\in p(\mathfrak M) \ |\ V^\psi_{p(\mathfrak M)}(\alpha)<\gamma<V^\psi_{p(\mathfrak M)}(\beta)\}$.
\end{center}
It is obvious that $\varepsilon_{\varphi_{min},\Theta}(\alpha)\subset V^\psi_{p(\mathfrak M)}(\alpha)\subset
\varepsilon_{\psi,\Theta}(\alpha)\subset\varepsilon_{\varphi_{max},\Theta}(\alpha)$.

If we treat $V^\psi_{p(\mathfrak M)}$ as $V_{p(\mathfrak M)}$ from the proof of Theorem~\ref{NoGreatest},
and work inside $\varepsilon_{\psi,\Theta}(\alpha)$ for some $\alpha\in p(\mathfrak M)$, we obtain an
analogical result to Lemma~\ref{MoreConvexTypes}.

\begin{lemma}\label{MoreConvexTypesPsi}
Let $n<\omega$, $\alpha\in p(\mathfrak M)$ and $\alpha_i\in \varepsilon_{\psi,\Theta}(\alpha)$, $i\leq n$,
such that $V^\psi_{p(\mathfrak M)}(\alpha_i)<V^\psi_{p(\mathfrak M)}(\alpha_{i+1})$, $1\leq i\leq n-1$.
If $\gamma_1,\gamma_2\in \varepsilon_{\psi,\Theta}(\alpha)$ are such that 
for some $i$, $1\leq i\leq n-1$,
$\gamma_1,\gamma_2\in (V^\psi_{p(\mathfrak M)}(\alpha_i),V^\psi_{p(\mathfrak M)}(\alpha_{i+1}))_p$, then
$tp^c(\gamma_1/A \cup\{\alpha_1,...,\alpha_n\})=
tp^c(\gamma_2/A\cup\{\alpha_1,...,\alpha_n\})$.
\end{lemma}

Lemma~\ref{MoreConvexTypesPsi} implies that the type $tp(\gamma_1/A\cup\{\alpha_1,...,\alpha_n\})$ is non-principal.

Let $\mathfrak N$ be an $\aleph_1$-saturated elementary extension of $\mathfrak M$.
Fix $\alpha\in p(\mathfrak N)$.
Let $\tau:=\langle\tau_1, \tau_2, ..., \tau_i,...\rangle_{i<\omega}$,
$\tau_i\in\{0,1\}$, be an infinite sequence of zeroes and ones.
Let $B_\tau:=\{\beta_1,\beta_2\}\cup\{\beta_{2n-1,q}\ |\ n\in\mathbb N,q\in \mathbb Q\}\cup
\{\beta_{2n,1}, \beta_{2n,2} \ |\ n\in \mathbb N, \tau_n=0\}\cup$
$\{\beta_{2n,1}, \beta_{2n,2}, \beta_{2n,3}\ |\ n\in \mathbb N, \tau_n=1\}$
be a subset of $\varepsilon_{\psi,\Theta}(\alpha)$ such that the sets $V^\psi_{p(\mathfrak N)}(\beta_{n,q})$
are disjoint and ordered lexicographically by the indices $n,q$,
and $V^\psi_{p(\mathfrak N)}(\beta_1)<V^\psi_{p(\mathfrak N)}(\beta_{n,q})<
V^\psi_{p(\mathfrak N)}(\beta_2)$ for every $\beta_{n,q}$.

Apply Theorem~\ref{Construction} to construct a model $\mathfrak M_\tau:=\mathfrak A^{B_\tau}$.

We claim that
$\varepsilon_{\psi,\Theta}(M_\tau,\alpha)\setminus
\underset{\beta\in B\tau}{\bigcup}V^\psi_{p(\mathfrak M_\tau)}(\beta)=\emptyset$.
As in Theorem~\ref{NoGreatest}, towards a contradiction, suppose that there exists a realization
$\gamma\in \varepsilon_{\psi,\Theta}(M_\tau,\alpha)\setminus\underset{\beta\in B_\tau}{\bigcup}V^\psi_{p(\mathfrak M_\tau)}(\beta)$.
By Theorem~\ref{Construction}, 3), the type $tp(\gamma/\bar \beta')$ is
principal for some $\bar \beta'\in B_\tau$ containing $\beta_1$ and $\beta_2$.
But by Lemma~\ref{MoreConvexTypesPsi}
the type $tp(\gamma/\bar \beta')$ is non-principal.
This is a contradiction and we have that
$p(\mathfrak M_\tau)\setminus\underset{\beta\in B_\tau}{\bigcup}V_{p(\mathfrak M_\tau)}(\beta)=\emptyset$.
Then we have as many countable non-isomorphic models of $I(T\cup tp(\bar a,\alpha),\aleph_0)$,
as there are different sequences $\tau$. Therefore, $I(T,\aleph_0)=2^{\aleph_0}$.

{\bf Case 2.}
There is a formula $\psi(x,y)\in CRF(p)$ such that for every $\psi_1\in CRF(p)$ with
$\psi(M,\alpha)\cap p(\mathfrak M)\subset \psi_1(M,\alpha)\cap p(\mathfrak M)$
for some (equivalently, for every) $\alpha\in p(\mathfrak M)$
there is $\psi_2\in CRF(p)$ such that  $\psi(M,\alpha)\subset\psi_2(M,\alpha) \subset\psi_1(M,\alpha)$.

For $\alpha\in p(\mathfrak M)$ denote $Ker^\psi_{p(\mathfrak M)}(\alpha):=
\underset{\varphi(M,\alpha)\cap p(\mathfrak M)\supset\psi(M,\alpha)\cap p(\mathfrak M)}
{\underset{\varphi\in CRF(p),}{\bigcap}}\varepsilon_{\varphi,\Theta}(M,\alpha)$.
It is obvious that $\varepsilon_{\varphi_{min},\Theta}(\alpha)
\subset\varepsilon_{\psi,\Theta}(\alpha)\subset
Ker^\psi_{p(\mathfrak M)}(\alpha)\subset\varepsilon_{\varphi_{max},\Theta}(\alpha)$.

We treat $Ker^\psi_{p(\mathfrak M)}$ as $Ker_{p(\mathfrak M)}$ from the proof of Theorem~\ref{NoLeast}
and work inside $\varepsilon_{\varphi_{max},\Theta}(\alpha)$ for a some $\alpha\in p(\mathfrak M)$.

\begin{lemma}\label{OneKernelPsi}
Let $\alpha,\alpha_1,\alpha_2\in p(\mathfrak M)$, $\alpha_1\neq\alpha_2\in Ker^\psi_{p(\mathfrak M)}(\alpha_1)\subset
\varepsilon_{\varphi_{max},\Theta}(M,\alpha)$,
and $\alpha_2\not\in\varepsilon_{\psi,\Theta}(M,\alpha_1)$.
Then $tp(\alpha_2/\alpha_1)$ is non-principal.
\end{lemma}

\begin{lemma}\label{ManyKernelsPsi}
Let $\alpha_1,...,\alpha_n, \alpha_{n+1}, \beta_1,...,\beta_m\in \big(p(\mathfrak M)\cap
\varepsilon_{\varphi_{max},\Theta}(M,\alpha)\big)$ ($n,m<\omega$)
be elements from different $\varepsilon_{\psi,\Theta}$-classes,
$\alpha_i\in Ker^\psi_{p(\mathfrak M)}(\alpha_{n+1})$ ($i\leq n$),
 and $\beta_j\not\in Ker^\psi_{p(\mathfrak M)}(\alpha_{n+1})$ ($j\leq m$).
Then $tp(\alpha_{n+1}/\bar\alpha\bar\beta)$ is non-principal, where
$\bar\alpha:=(\alpha_1,...,\alpha_n)$, $\bar\beta:=(\beta_1,...,\beta_m)$.
\end{lemma}

Fix $\alpha\in p(\mathfrak N)$.
Let $\tau:=\langle\tau_1, \tau_2, ..., \tau_i,...\rangle_{i<\omega}$,
$\tau_i\in\{0,1\}$, be an infinite sequence of zeroes and ones.
Let $B_\tau:=\{\beta_{2n-1,q,1}\ |\ n\in\mathbb N, q\in\mathbb Q\}\cup
\{\beta_{2n,1,m}\ |\ n,m\in\mathbb N, \tau_n=0\}\cup
\{\beta_{2n,1,m}\ |\ n\in\mathbb N, m\in\mathbb N^- \tau_n=1\}$
be a subset of $p(\mathfrak N)$ such that all $\beta_{n_1,q_1,m_1}$ are in $\varepsilon_{\varphi_{max},\Theta}(\alpha)$;
$\beta_{n_1,q_1,m_1}\in Ker^\psi_{p(\mathfrak N)}(\beta_{n_2,q_2,m_2})$
if and only if $n_1=n_2$ and $q_1=q_2$;
all the $\beta_{n,q,m}$ are in different $\varepsilon_{\psi,\Theta}$-classes and
 are ordered lexicographically by their indices;
and for every element from $p(\mathfrak N)\cap \varepsilon_{\varphi_{max},\Theta}(\alpha)$
there is a representative of its kernel in $B_\tau$.
Apply Theorem~\ref{Construction} to construct a model $\mathfrak M_\tau:=\mathfrak A^{B_\tau}$.

Suppose  there exists $\gamma\in \varepsilon_{\varphi_{max},\Theta}(N,\alpha)\setminus B_\tau$
such that for some $\alpha'\in B_\tau$ $\gamma\in Ker^\psi_{p(\mathfrak (M_\tau)}(\alpha')$
and $\gamma>\varepsilon_{\psi,\Theta}(M_\tau,\alpha'')$ ($\gamma<\varepsilon_{\psi,\Theta}(M_\tau,\alpha'')$)
for every $\alpha''\in B_\tau\cap Ker^\psi_{p(\mathfrak (M_\tau)}(\alpha')$.
By Theorem~\ref{Construction}, 3), there is $\bar \beta\in B_\tau$ such that $tp(\gamma/\bar \beta)$ is principal.
By the proofs of Lemma~\ref{chain} and Theorem~\ref{Construction}
there is $\bar \beta'\in B_\tau$ extending $\bar \beta$ and
such that $\bar \beta'$ includes at least one element form $Ker_{p(\mathfrak M_\tau)}(\alpha')\cap B_\tau$
and such that $tp(\gamma/\bar \beta')$ is also principal.
But by Lemma~\ref{ManyKernelsPsi} $tp(\gamma/\bar b')$ is non-principal.
This is a contradiction, and, therefore, there is no such $\gamma$.

Analogically to Theorem~\ref{NoLeast}, by Lemma~\ref{MoreConvexTypesPsi},
 the order types of kernels in $\varepsilon_{\varphi_{max},\Theta}(N,\alpha)$
are determined by their order types in $B_\tau$, and
$\mathfrak M_{\tau_1}\not\cong\mathfrak M_{\tau_2}$ whenever $\tau_1\neq\tau_2$.

Since the number of different sequences $\tau$ equals to
$2^{\aleph_0}$, $I(T\cup tp(\bar a,\alpha),\aleph_0)=2^{\aleph_0}$,
and, therefore, $I(T,\aleph_0)=2^{\aleph_0}$.

{\bf Case 3.} Suppose that the cases 1 and 2 are not true. Then there exists a partition 
$CRF(p)=L\cup R$ such that
for some (equivalently, for every) $\alpha\in p(\mathfrak M)$, for every $\varphi_l\in L$ and every $\varphi_r\in R$,
$\varphi_l(M,\alpha)\cap p(\mathfrak M)\subset\varphi_r
(M,\alpha)\cap p(\mathfrak M)$.

Denote $Ker_{p(\mathfrak M)}(\alpha):=
\underset{\varphi\in R}{\bigcap}\varepsilon_{\varphi,\Theta}(M,\alpha)$;
$V_{p(\mathfrak M)}(\alpha):=\underset{\varphi\in L}{\bigcup}
\varepsilon_{\varphi,\Theta}(M,\alpha)$.
Then $V_{p(\mathfrak M)}(\alpha)\subset Ker_{p(\mathfrak M)}(\alpha)$.
For $\tau\in 2^\omega$ let
$B_\tau:=\{\beta_{2n-1,q}\ |\ n\in\mathbb N, q\in Q\}\cup \{\beta_{2n,1},\beta_{2n,2}\ |\ n\in\mathbb N, \tau_n=0\}\cup
\{\beta_{2n,1},\beta_{2n,2},\beta_{2n,3}\ |\ n\in\mathbb N, \tau_n=1\}$ be a subset of 
$p(\mathfrak N)$ such that all the $\beta_{q,n}$'s are in different $V_{p(\mathfrak M)}$-neighborhoods,
lexicographically ordered by the indices, and $\beta_{n_1,q_1}\in Ker(\beta_{n_2,q_2})$ if and only if $n_1=n_2$.
Apply Theorem~\ref{Construction} to construct a model $\mathfrak M_\tau:=\mathfrak A^{B_\tau}$.
As in the previous cases, the models constructed for different sequences $\tau$ are pairwise non-isomorphic
and therefore the theorem is proved.

$\hfill\Box$ Theorem~\ref{Infinite}

\noindent{}

\begin{footnotesize}
\textsc{Suleyman Demirel University,}

\textsc{Kaskelen, Kazakhstan}

\textsc{Institute of Mathematics and Mathematical Modeling,}

\textsc{Almaty, Kazakhstan}

{\it Email address:} baizhanov@math.kz

\textsc{Institute of Mathematics and Mathematical Modeling,}

\textsc{Almaty, Kazakhstan}

{\it Email address:} zambarnaya@math.kz
\end{footnotesize}

\end{document}